\documentclass[10pt]{article}
\usepackage{amsthm, amsmath, amsfonts,	xcolor,	hyperref,	enumerate,	geometry, multirow,	caption, orcidlink
}
\numberwithin{equation}{section}
\newtheorem{theorem}{Theorem}[section]
\newtheorem{example}[theorem]{Example}

\title{\textbf{General Chen-Ricci inequalities for Riemannian submersions and Riemannian maps}}

\author{Ravindra Singh\orcidlink{0009-0009-1270-3831}
	, Kiran Meena\orcidlink{0000-0002-6959-5853}
	, Kapish Chand Meena\orcidlink{0000-0003-0182-8822}
}

\date{}

\begin{document}
	\maketitle
	\begin{abstract}
		\noindent In this paper, we derive general forms of the Chen-Ricci inequalities for Riemannian submersions between Riemannian manifolds. We also derive general forms of the Chen-Ricci and improved Chen-Ricci inequalities for Riemannian maps between Riemannian manifolds, involving relations between the curvatures of subspaces of the source and target spaces. Further, we illustrate equality cases for all these general forms with two examples. These general forms yield new, easy, and elegant techniques that are fruitful in obtaining the Chen-Ricci inequalities for such smooth mappings with various structured manifolds. As applications, utilizing these general forms, we explicitly establish Chen-Ricci inequalities when the source manifolds of Riemannian submersions and the target manifolds of Riemannian maps belong to broader classes, such as generalized complex and generalized Sasakian space forms, particularly including real, complex, real K\"ahler, Sasakian, Kenmotsu, cosymplectic, and almost $C(\alpha)$ space forms. We also validate our approach by imposing appropriate conditions toward various particular existing cases.
	\end{abstract}
	\noindent \textbf{Keywords: }{Riemannian submersions; Riemannian maps; (improved) Chen-Ricci inequalities; complex and contact manifolds; space forms.}\\
	
	\noindent \textbf{MSC Classification:} {53B20, 53B35, 53C15, 53D15}\\
	
	
	\section{Introduction}

    In submanifold or immersion theory, one of the most fundamental problems is to establish simple and sharp relationships between the intrinsic and extrinsic curvature invariants of a Riemannian submanifold \cite{Chen_2011}. Toward this, various authors have derived \textit{Chen-Ricci inequalities} providing relations between the Ricci curvature (intrinsic invariant) and the squared mean curvature (extrinsic invariant) for different types of submanifolds in several space forms. These inequalities were also improved by Oprea \cite{Oprea_2010} and Deng \cite{Deng_2009}, using optimization techniques and some algebraic inequalities, respectively. However, improved Chen-Ricci inequalities for K\"ahlerian slant submanifolds of complex space forms and for $C$-totally real submanifolds of Sasakian space forms cannot be obtained using the approaches of Oprea and Deng. To avoid such limitations, Tripathi \cite{Tripathi_2011} established an approach to improve Chen-Ricci inequalities. For a detailed survey on Chen-Ricci inequalities, we refer to \cite{Chen_Blaga} and the references therein.
	
	We know that a Riemannian submersion is the dual concept of isometric immersion and has various applications \cite{Falcitelli_2004, Sahin_book}. The concept of a Riemannian map encompasses the concepts of isometric immersion and Riemannian submersion, providing frameworks for quantum models and for comparing the geometric properties of two arbitrary Riemannian manifolds \cite{Fischer_1992, Sahin_book}. Since Riemannian submersions and Riemannian maps help to establish relations between curvature invariants of submanifolds of the source and target manifolds, the Chen-Ricci inequalities have also been studied by many authors in the context of these smooth mappings. For geometrically structured Riemannian submersions, they are explored from various space forms, as listed in Table \ref{table_rs}. In addition, for Riemannian maps with real and complex space forms as the target spaces, Lee et al. \cite{LLSV_2022} obtained Chen-Ricci and improved Chen-Ricci inequalities (using Deng's approach).
	\begin{table}[h]
		\centering
		\caption{Chen-Ricci inequalities for Riemannian submersions}\label{table_rs}
		\begin{tabular}{lll} 
			\hline
			\textbf{Space Form} & \textbf{Type of Submersion} & \textbf{Ref.}\\
			\hline			
			\multirow{4}{*}{complex} & Lagrangian & \cite{Gulbahar_Meric_Kilic}\\
			&anti-invariant & \cite{Gunduzalp_Polat_MMN}\\
			&slant & \cite{Gunduzalp_Polat_Filomat}\\
			&hemi-slant & \cite{Akyol_Demir_Poyraz_Vilcu}\\
			\hline
			\multirow{2}{*}{Sasakian} &	anti-invariant & \cite{Aytimur_Ozgur}\\
			&slant & \cite{Akyol_Poyraz}\\
			\hline
			\multirow{2}{*}{Kenmotsu} &anti-invariant & \cite{Polat_2024}\\
			&slant & \cite{Aquib_Aldayel_Iqbal_Khan}\\
			\hline
			\multirow{2}{*}{cosymplectic} & anti-invariant & \cite{Aytimur_2023}\\
			&slant & \cite{Poyraz_Akyol}\\
			\hline
			slant & generalized Sasakian & \cite{Aquib_2025}\\
			\hline 
		\end{tabular}
	\end{table} 
	
	Motivatingly, in the current paper, we introduce a general framework of Chen-Ricci inequalities for Riemannian submersions between Riemannian manifolds (Theorem $\ref{Theorem-GCRV}$). We also introduce general frameworks of the Chen-Ricci and improved Chen-Ricci inequalities (using Tripathi's approach) for Riemannian maps between Riemannian manifolds (Theorems  $\ref{main_thm_CRI}$ and  $\ref{main_thm_ICRI}$). Notably, these frameworks unify existing studies and yield relations between the curvatures of subspaces of source and target spaces. These frameworks are useful for deriving Chen-Ricci inequalities for Riemannian submersions from arbitrary source manifolds and for Riemannian maps to arbitrary target manifolds (for instance, see Theorems $\ref{appl_for_CRI_RS}$, $\ref{appl_thm_CRI}$ and $\ref{appl_thm_ICRI}$). In addition, we can use these general frameworks to establish Chen-Ricci inequalities for various geometrically structured Riemannian submersions and Riemannian maps, such as invariant, anti-invariant, Lagrangian, semi-invariant, slant, semi-slant, hemi-slant, etc. (see Section \ref{sec_conclusion}).
	
	\section{Generalized space forms}\label{sec_space_forms}
	In this section, we recall the notions of generalized space forms with their special cases.
	
	\subsubsection*{Generalized complex space forms}
	Let $({M}, g, J)$ be an even-dimensional almost Hermitian manifold. It is referred to as a {\it generalized complex space form}, denoted by ${M}(c_{1},c_{2})$, if its curvature tensor $R^{M}$ satisfies $\cite{Olszak_1989, Tricerri_1981, Yano_1984}$, 
	\begin{align}\label{curvature_for_gcsf}
		R^{M}({\cal Y}_{1},{\cal Y}_{2}){\cal Y}_{3}& =c_{1}\{g({\cal Y}_{2},{\cal Y}_{3}){\cal Y}_{1}-g({\cal Y}_{1},{\cal Y}_{3}){\cal Y}_{2}\} \\& \quad +c_{2}\{g({\cal Y}_{1},J{\cal Y}_{3}) J{\cal Y}_{2}-g({\cal Y}_{2},J{\cal Y}_{3})J{\cal Y}_{1}+2g({\cal Y}_{1},J{\cal Y}_{2})J{\cal Y}_{3}\}\nonumber
	\end{align}
	for all ${\cal Y}_{1}, {\cal Y}_{2}, {\cal Y}_{3} \in \Gamma (T{M})$. Here, $c_{1}$ and $c_{2}$ are smooth functions on ${M}$. 
	
	\subsubsection*{Generalized Sasakian space forms}
	Let $({M}, g, \phi, \xi, \eta)$ be an odd-dimensional almost contact metric manifold. It is termed as {\it generalized Sasakian space form} (denoted by $M(c_{1},c_{2},c_{3})$) if there exist smooth functions $c_{1}$, $c_{2}$, and $c_{3}$ on $M$ such that the curvature tensor $R^{M}$ of $M$ satisfies $\cite{Alegre_Blair, Blair_2010}$, 
	\begin{align}\label{curvature_for_gssf}
		&R^{M}({\cal Y}_{1},{\cal Y}_{2}){\cal Y}_{3} \nonumber \\&= c_{1}\{g({\cal Y}_{2},{\cal Y}_{3}){\cal Y}_{1}-g({\cal Y}_{1},{\cal Y}_{3}){\cal Y}_{2}\} \nonumber \\& \quad +c_{2}\{g({\cal Y}_{1},\phi {\cal Y}_{3})\phi {\cal Y}_{2}-g({\cal Y}_{2},\phi {\cal Y}_{3})\phi {\cal Y}_{1}+2g({\cal Y}_{1},\phi {\cal Y}_{2})\phi {\cal Y}_{3}\} \\& \quad +c_{3}\{\eta ({\cal Y}_{1})\eta ({\cal Y}_{3}){\cal Y}_{2}-\eta ({\cal Y}_{2})\eta ({\cal Y}_{3}){\cal Y}_{1}+g({\cal Y}_{1},{\cal Y}_{3})\eta ({\cal Y}_{2})\xi -g({\cal Y}_{2},{\cal Y}_{3})\eta ({\cal Y}_{1})\xi \},\nonumber
	\end{align}
	for any ${\cal Y}_{1}, {\cal Y}_{2}, {\cal Y}_{3} \in \Gamma (T{M})$.\\
	
	By \cite{Alegre_Blair, Blair_2010, Vanhecke_1975, Yano_1984, Tricerri_1981}, these generalized space forms include several special space forms based on the particular values of functions $c_1$, $c_2$, and $c_3$ as follows.
	
	\begin{table}[!h]
		\centering
		\caption{Special cases}\label{table_3}
		\begin{tabular}{cccc}
			\hline
			\textbf{space form}& \textbf{$c_1$} & \textbf{$c_2$} & \textbf{$c_3$} \\
			\hline 
			real & $c$ & $0$ & -\\
			\vspace{0.1cm}
			
			complex & $\frac{c}{4}$ & $\frac{c}{4}$ & -\\
			\vspace{0.1cm}
			real K\"ahler &$\frac{c + 3\alpha}{4}$& $\frac{c - \alpha}{4}$ & -\\
			\vspace{0.1cm} 
			Sasakian & $\frac{c + 3}{4}$ & $\frac{c - 1}{4}$& $\frac{c - 1}{4}$\\
			\vspace{0.1cm}
			Kenmotsu & $\frac{c - 3}{4}$ & $\frac{c + 1}{4}$& $\frac{c + 1}{4}$\\
			\vspace{0.1cm}
			cosymplectic & $\frac{c}{4}$ & $\frac{c}{4}$& $\frac{c}{4}$\\
			\vspace{0.1cm}
			almost $C(\alpha)$ &$\frac{c + 3\alpha^2}{4}$& $\frac{c-\alpha^2}{4}$& $\frac{c-\alpha^2}{4}$\\
			\hline
		\end{tabular}
	\end{table}
	
	\section{General Chen-Ricci inequalities for Riemannian submersions with applications}\label{part_I}
	\subsection{Background}\label{sec_prilims}
	In this subsection, we recall some information required for the next subsections.\\
	
	Let $\left( M_{1}^{m_{1}},g_{1}\right)$ and $\left(M_{2}^{m_{2}},g_{2}\right) $ be two Riemannian manifolds. A smooth surjective map $F:\left( M_{1}^{m_{1}},g_{1}\right)$ $\to \left(M_{2}^{m_{2}},g_{2}\right)$ is called a {\it Riemannian submersion} if its differential map (i.e., $F_{\ast p}: T_{p} M_{1} \to T_{F(p)}M_{2}$) is surjective at all $p \in M_1$ and $F_{\ast p}$ preserves the length of all horizontal vectors at $p$. Here, the tangent vectors to the fibers (fiber means $\{F^{-1}(q)$ for $q\in M_2$\}) are termed {\it vertical}, while the orthogonal vectors to the fibers are termed {\it horizontal}. This induces a decomposition of the tangent bundle $T M_1$ into a direct sum of two distributions: the \textit{vertical distribution} $\mathcal{V} = \ker F_{\ast}$ and its orthogonal complement (known as the \textit{horizontal distribution}) $\mathcal{H} = \left( \ker F_{\ast} \right)^{\perp}$.
	\subsubsection*{O'Neill tensors}
	For a Riemannian submersion, O'Neill \cite{Neill_1996} defined two $(1,2)$-type tensor fields ${T}$ and ${A}$ that satisfy
	\begin{equation*}
		{A}_{Y_{1}}Y_{2}=-A_{Y_{2}}Y_{1}, 
	\end{equation*}
	\begin{equation*}
		{T}_{U_{1}}U_{2}={T}_{U_{2}}U_{1}, 
	\end{equation*}
	\begin{equation*}
		g_{1}\left( {T}_{U_{1}}Y_{2},Y_{3}\right) =-g_{1}\left( {T}_{U_{1}}Y_{3},Y_{2}\right),
	\end{equation*}
	and
	\begin{equation*}
		g_{1}\left( {A}_{Y_{1}}Y_{2},Y_{3}\right) =-g_{1}\left( {A}_{Y_{1}} Y_{3}, Y_{2}\right),
	\end{equation*}
	where $U_{1}, U_{2}\in \chi \left(\ker F_{\ast }\right)$ and $Y_{1}, Y_{2}, Y_{3} \in \chi (\left( \ker F_{\ast }\right)^{\perp })$. 
	Moreover, let $\left\{ \vee_{1},\ldots, \vee_{n}\right\} $ be an orthonormal basis of the vertical distribution $(\ker F_{\ast})$. Then the \textit{mean curvature vector field $H$ of the fibers} of $F$ is defined as \cite{Falcitelli_2004}, 
	\begin{equation}\label{mean_curv}
		H= \frac{1}{n} \sum \limits_{i=1}^{n}T_{\vee_{i}}\vee_{i}.
	\end{equation}
	\subsubsection*{Relations between Riemannian curvature tensors}
	Let $R^{M_{1}}$, $R^{M_{2}}$, $R^{\ker F_{\ast }}$, and $R^{\left( \ker F_{\ast }\right)^{\perp }}$ denote the Riemannian curvature tensors corresponding to $M_{1}$, $M_{2}$, $\ker F_{\ast }$, and $\left( \ker F_{\ast }\right)^{\perp }$, respectively. Then we have \cite{Neill_1996, Falcitelli_2004}
	\begin{eqnarray}\label{eq-(2.3)}
		R^{M_{1}}\left( U_{1},U_{2},U_{3},U_{4}\right) &=&R^{\ker F_{\ast }}\left(U_{1},U_{2},U_{3},U_{4}\right) +g_{1}\left( {T}_{U_{1}}U_{4},{T}_{U_{2}}U_{3}\right) \nonumber \\&&-g_{1}\left( {T}_{U_{2}}U_{4},{T}_{U_{1}}U_{3}\right),
	\end{eqnarray}
	\begin{eqnarray}\label{eq-(2.4)}
		R^{M_{1}}\left( Y_{1},Y_{2},Y_{3},Y_{4}\right) &=&R^{\left( \ker F_{\ast}\right)^{\perp }}\left( Y_{1},Y_{2},Y_{3},Y_{4}\right) -2g_{1}\left({A}_{Y_{1}}Y_{2},{A}_{Y_{3}}Y_{4}\right) \nonumber \\&&+g_{1}\left( {A}_{Y_{2}}Y_{3},{A}_{Y_{1}}Y_{4}\right)-g_{1}\left({A}_{Y_{1}}Y_{3},{A}_{Y_{2}}Y_{4}\right),
	\end{eqnarray}
	and
	\begin{eqnarray}\label{eq-(2.5)}
		R^{M_{1}}\left( Y_{1},U_{1},Y_{2},U_{2}\right) &=&g_{1}\left( \left( \nabla_{Y_{1}}^{1}{T}\right) \left( U_{1},U_{2}\right), Y_{2}\right)+g_{1}\left( \left( \nabla_{U_{1}}^{1}{A}\right) \left(Y_{1},Y_{2}\right), U_{2}\right) \nonumber \\&&-g_{1}\left({T}_{U_{1}}Y_{1},{T}_{U_{2}}Y_{2}\right)+g_{1}\left( {A}_{Y_{2}}U_{2},{A}_{Y_{1}}U_{1}\right),
	\end{eqnarray}
	for all $Y_{1},Y_{2},Y_{3},Y_{4} \in \chi \left( \left( \ker F_{\ast }\right)^{\perp }\right) $ and $U_{1},U_{2},U_{3},U_{4}\in \chi \left( \ker F_{\ast}\right)$. Here, $\nabla^1$ is the Levi-Civita connection with respect to the metric $g_1$.
	
	\subsubsection*{Some notations}
	Let $(\ker F_{\ast }) = {\rm span}\{\vee_{1},\ldots, \vee_{n}\}$ and $\left(\ker F_{\ast }\right)^{\perp }={\rm span}\{h_{1},\ldots, h_{r}\}$. We define some notation as
	\begin{equation}\label{eq-(3.1)}
		{T}_{ij}^{t} := g_1({T}_{\vee_i} \vee_j, h_t), \quad 1 \leq i,j \leq n,\
		1 \leq t \leq r, 
	\end{equation}
	\begin{equation}\label{eq-(3.2)}
		A_{ij}^{\alpha} := g_1({A}_{h_i} h_j, \vee_{\alpha}), \quad 1 \leq i,j
		\leq r,\ 1 \leq \alpha \leq n,
	\end{equation}
	\begin{equation}\label{eq-(3.3)}
		\delta(N) := \sum \limits_{i=1}^{r} \sum \limits_{j=1}^{n} \left( \left( \nabla^{1}_{h_i} 
		{T} \right)_{\vee_j} \vee_j, h_i \right),
	\end{equation}
	\begin{equation}\label{TVert}
		\left\Vert T^{\cal V}\right\Vert^{2}:= \sum \limits_{i=1}^{r} \sum \limits_{j=1}^{n}g_{1}\left(T_{\vee_{j}}h_{i},T_{\vee_{j}}h_{i}\right),
	\end{equation}
	and
	\begin{equation}\label{AHor}
		\left\Vert A^{{\cal H}}\right\Vert^{2}:= \sum \limits_{i=1}^{r} \sum \limits_{j=1}^{n}g_{1}\left(A_{h_{i}}\vee_{j},A_{h_{i}}\vee_{j}\right).
	\end{equation}
	Then from \cite{Gulbahar_Meric_Kilic}, we have 
	\begin{align}\label{eq-(3.5)}
		\sum \limits_{t=1}^{r} \sum \limits_{i,j=1}^{n}\left( T_{ij}^{t}\right)^{2}& =\frac{1}{2} n^{2}\left\Vert H\right\Vert^{2}+\frac{1}{2} \sum \limits_{t=1}^{r}\left(T_{11}^{t}-T_{22}^{t}-\cdots -T_{nn}^{t}\right)^{2} \nonumber \\& +2 \sum \limits_{t=1}^{r} \sum \limits_{j=2}^{n}\left( T_{1j}^{t}\right)^{2}-2 \sum \limits_{t=1}^{r} \sum \limits_{2\leq i<j\leq n}\left\{T_{ii}^{t}T_{jj}^{t}-\left( T_{ij}^{t}\right)^{2}\right\}. 
	\end{align}
	In the sequel, we also use curvature-like notation. Therefore, we fix those as follows.
	\begin{equation*}
		2\tau_{{\cal H}}^{\left( \ker F_{\ast }\right)^{\perp }} = \sum \limits_{i,j=1}^{r}R^{\left( \ker F_{\ast }\right)^{\perp }}\left(h_{i},h_{j},h_{j},h_{i}\right), \quad 2\tau_{{\cal H}}^{M_{1}} = \sum \limits_{i,j=1}^{r}R^{M_{1}}\left(h_{i},h_{j},h_{j},h_{i}\right),
	\end{equation*}
	\begin{equation*}
		2\tau_{{\cal V}}^{\ker F_{\ast }} = \sum \limits_{i,j=1}^{n}R^{\ker F_{\ast }}\left(\vee_{i},\vee_{j},\vee_{j},\vee_{i}\right), \quad 2\tau_{{\cal V}}^{M_{1}} = \sum \limits_{i,j=1}^{n}R^{M_{1}}\left(\vee_{i},\vee_{j},\vee_{j},\vee_{i}\right),
	\end{equation*}
	\begin{equation*}
		{\rm Ric}_{{\cal H}}^{\left( \ker F_{\ast }\right)^{\perp }} (h_1) = \sum \limits_{j=1}^{r}R^{\left( \ker F_{\ast }\right)^{\perp }}\left(h_{1},h_{j},h_{j},h_{1}\right), \quad {\rm Ric}_{{\cal H}}^{M_{1}} (h_1) = \sum \limits_{j=1}^{r}R^{M_{1}}\left(h_{1},h_{j},h_{j},h_{1}\right),
	\end{equation*}
	\begin{equation*}
		{\rm Ric}_{{\cal V}}^{\ker F_{\ast }} (\vee_1) = \sum \limits_{j=1}^{n}R^{\ker F_{\ast }}\left(\vee_{1},\vee_{j},\vee_{j},\vee_{1}\right), \quad {\rm Ric}_{{\cal V}}^{M_{1}}(\vee_1) = \sum \limits_{j=1}^{n}R^{M_{1}}\left(\vee_{1},\vee_{j},\vee_{j},\vee_{1}\right).
	\end{equation*}
	
	\subsection{General Chen-Ricci inequalities}\label{sec_general_CRI}
	Here, we present general forms of the Chen-Ricci inequalities for Riemannian submersions between Riemannian manifolds. Toward this, we have the following result. 
	
	\begin{theorem}\label{Theorem-GCRV}
		Let $F:(M_{1}^{m_{1}},g_{1}) \to (M_{2}^{m_{2}},g_{2})$ be a Riemannian submersion between two Riemannian manifolds. Then
		\begin{equation}\label{eq-GFCRV-(1)}
			{\rm Ric}_{{\cal V}}^{\ker F_{\ast }}(\vee_{1})\geq {\rm Ric}_{{\cal V}}^{M_{1}}(\vee_{1})-\frac{1}{4}n^{2}\left\Vert H\right\Vert^{2},
		\end{equation}
		
		\begin{equation}\label{eq-GFCRH-(1)}
			{\rm Ric}_{{\cal H}}^{\left( \ker F_{\ast }\right)^{\perp}}(h_{1})\leq {\rm Ric}_{{\cal H}}^{M_{1}}(h_{1}), 
		\end{equation}
		and
		\begin{align}\label{eq-GFCRVH-(1)}
			{\rm Ric}_{{\cal V}}^{\ker F_\ast}\left( \vee_{1}\right) & + {\rm Ric}_{\cal H}^{(\ker F_\ast)^\perp}\left(h_{1}\right) \nonumber \\& \geq {\rm Ric}_{{\cal V}}^{M_{1}} \left( \vee_{1}\right) +{\rm Ric}_{{\cal H}}^{M_{1}} \left( h_{1}\right) + \sum \limits_{i=1}^{r} \sum \limits_{j=1}^{n}R^{M_{1}}\left(h_{i},\vee_{j},\vee_{j},h_{i}\right) \nonumber \\&+\delta \left( N\right) - \left\Vert T^{{\cal V}}\right\Vert^{2} + \left\Vert A^{{\cal H}}\right\Vert^{2} -\frac{n^{2}}{4}\left\Vert H\right\Vert^{2} - 3 \sum \limits_{\alpha=1}^{n} \sum \limits_{j=2}^{r}\left( A_{1j}^{\alpha }\right)^{2}.
		\end{align}
		Further, equality cases are as follows. \\
		\begin{enumerate}[$(1)$]
			\item 
			The equality in $(\ref{eq-GFCRV-(1)})$ and $(\ref{eq-GFCRVH-(1)})$ is achieved if and only if 
			\[
			T_{11}^{t}=T_{22}^{t}+\cdots +T_{nn}^{t} \quad \text{and}\quad T_{1j}^{t}=0 \quad \text{for} \quad j \in \{2,\ldots, n\},\quad t \in \{1,\ldots, r\}.
			\]
			Equivalently, the equality in $(\ref{eq-GFCRV-(1)})$ and $(\ref{eq-GFCRVH-(1)})$ is achieved for all unit vertical vectors if and only if the fibers of $F$ are
			\begin{enumerate}[$(a)$]
				\item totally geodesic when $n>2$;
				
				\item totally umbilical when $n=2$.\\
			\end{enumerate}
			
			\item 
			The equality in $(\ref{eq-GFCRH-(1)})$ is achieved if and only if
			\begin{equation*}
				A_{1j}^{\alpha }=0 \quad \text{for} \quad j \in \{2,\ldots, r\},\quad \alpha \in \{1,\ldots, n\}.
			\end{equation*}
			Equivalently, the equality in $(\ref{eq-GFCRH-(1)})$ is achieved for all unit horizontal vectors if and only if the horizontal distribution is integrable.
		\end{enumerate}
		
	\end{theorem}
	
	\begin{proof}
		For any point $p\in M_{1}$, choose an orthonormal basis of $T_{p}M_{1}$ as $\{h_{1},\ldots, h_{r},\vee_{1},\ldots, \vee_{n}\}$, where $(\ker F_{\ast }) = {\rm span}\{\vee_{1},\ldots, \vee_{n}\}$ and $\left( \ker F_{\ast }\right)^{\perp } = {\rm span}\{h_{1}, \ldots, h_{r}\}$. Then from (\ref{eq-(2.3)}) and (\ref{mean_curv}), we obtain
		\begin{equation*}
			2\tau_{{\cal V}}^{\ker F_{\ast }}=2\tau_{{\cal V}}^{M_{1}}-n^{2}\left\Vert H\right\Vert^{2}+ \sum \limits_{i,j=1}^{n}g_{1}\left(T_{\vee_{j}}\vee_{i},T_{\vee_{i}}\vee_{j}\right).
		\end{equation*}
		Employing (\ref{eq-(3.1)}) and then substituting (\ref{eq-(3.5)}) yield
		\begin{align*}
			2\tau_{{\cal V}}^{\ker F_{\ast }}& =2\tau_{{\cal V}}^{M_{1}}-\frac{1}{2}n^{2}\left\Vert H\right\Vert^{2}+\frac{1}{2} \sum \limits_{t=1}^{r}\left( T_{11}^{t}-T_{22}^{t}-\cdots -T_{nn}^{t}\right)^{2}	\\& \quad +2 \sum \limits_{t=1}^{r} \sum \limits_{j=2}^{n}\left( T_{1j}^{t}\right)^{2}-2 \sum \limits_{t=1}^{r} \sum \limits_{2\leq i<j\leq n}\left\{T_{ii}^{t}T_{jj}^{t}-\left( T_{ij}^{t}\right)^{2}\right\}.
		\end{align*}
		Thus, 
		\begin{equation}\label{eq-(3.6)}
			2\tau_{{\cal V}}^{\ker F_{\ast }}\geq 2\tau_{{\cal V}}^{M_{1}}-\frac{1}{2} n^{2}\left\Vert H\right\Vert^{2}-2 \sum \limits_{t=1}^{r} \sum \limits_{2\leq i<j\leq n}\left\{ T_{ii}^{t}T_{jj}^{t}-\left( T_{ij}^{t}\right)^{2}\right\}.
		\end{equation}
		Using (\ref{eq-(2.3)}) and (\ref{eq-(3.1)}), we derive
		\begin{align}\label{eq-(3.6a)}
			2 \sum \limits_{2\leq i<j\leq n}R^{M_{1}}(\vee_{i},\vee_{j},\vee_{j},\vee_{i})=&2 \sum \limits_{2\leq i<j\leq n}R^{\ker F_{\ast}}(\vee_{i},\vee_{j},\vee_{j},\vee_{i}) \nonumber \\&+2 \sum \limits_{t=1}^{r} \sum \limits_{2\leq i<j\leq n}\left({T}_{ii}^{t}{T}_{jj}^{t}-({T}_{ij}^{t})^{2}\right).
		\end{align}
		Substituting (\ref{eq-(3.6a)}) into (\ref{eq-(3.6)}), we get 
		\begin{align}\label{eq-(3.7)}
			2\tau_{{\cal V}}^{\ker F_{\ast }}&-2 \sum \limits_{2\leq i<j\leq n} R^{\ker F_{\ast }}(\vee_{i},\vee_{j},\vee_{j},\vee_{i}) \nonumber \\& \geq 2\tau_{{\cal V}}^{M_{1}}-2 \sum \limits_{2	\leq i<j\leq n}R^{M_{1}}(\vee_{i},\vee_{j},\vee_{j},\vee_{i})-\frac{1}{2}n^{2}\Vert H \Vert^{2}.
		\end{align}
		Since $2\tau_{{\cal V}}^{\ker F_{\ast }}=2 \sum \limits_{2\leq i<j\leq n} R^{\ker F_{\ast}}(\vee_{i},\vee_{j},\vee_{j},\vee_{i}) + 2 \sum \limits_{j=1}^{n}R^{\ker F_{\ast}}(\vee_{1},\vee_{j},\vee_{j},\vee_{1})$, the equation (\ref{eq-(3.7)}) implies (\ref{eq-GFCRV-(1)}). In addition, its equality case is straightforward by (\ref{eq-(3.6)}).\\
		
		Further, utilizing (\ref{eq-(2.4)}) and (\ref{eq-(3.2)}), we obtain 
		\[
		\sum \limits_{i,j=1}^{r}R^{\left( \ker F_{\ast }\right)^{\perp }}\left(h_{i},h_{j},h_{j},h_{i}\right) = \sum \limits_{i,j=1}^{r}R^{M_{1}}\left(h_{i},h_{j},h_{j},h_{i}\right) -3 \sum \limits_{\alpha =1}^{n} \sum \limits_{i,j=1}^{r}\left(A_{ij}^{\alpha }\right)^{2}.
		\]
		Equivalently, 
		\begin{align}\label{eq-(4.2)}
			2\tau_{{\cal H}}^{\left( \ker F_{\ast }\right)^{\perp }} &= 2\tau_{{\cal H}}^{M_{1}} - 3 \sum \limits_{\alpha=1}^{n} \sum \limits_{i,j=1}^{r}\left( A_{ij}^{\alpha}\right)^{2} \nonumber \\& = 2 \tau_{{\cal H}}^{M_{1}} - 6 \sum \limits_{\alpha=1}^{n} \sum \limits_{j=2}^{r} \left(A_{1j}^{\alpha}\right)^{2} - 6 \sum \limits_{\alpha =1}^{n} \sum \limits_{2 \leq i<j\leq r} \left(A_{ij}^{\alpha}\right)^{2}.
		\end{align}
		Again, by substituting $Y_{1}=Y_{4}=h_{i}$ and $Y_{2}=Y_{3}=h_{j}$ in (\ref{eq-(2.4)}) and then using (\ref{eq-(3.2)}), we obtain 
		\begin{align}\label{eq-(4.3)}
			2 \sum \limits_{2\leq i<j\leq r}R^{M_{1}}\left( h_{i},h_{j},h_{j},h_{i}\right)=& 2 \sum \limits_{2 \leq i<j\leq r} R^{\left(\ker F_{\ast }\right)^{\perp }} \left(h_{i},h_{j},h_{j},h_{i}\right) \nonumber \\&+ 6 \sum \limits_{\alpha =1}^{n} \sum \limits_{2\leq i<j \leq r} \left(A_{ij}^{\alpha}\right)^{2}.
		\end{align}
		Substituting (\ref{eq-(4.3)}) into (\ref{eq-(4.2)}), we deduce that
		\begin{eqnarray*}
			2\tau_{{\cal H}}^{\left(\ker F_{\ast }\right)^{\perp}} & = & 2 \tau_{{\cal H}}^{M_{1}} - 6 \sum \limits_{\alpha=1}^{n} \sum \limits_{j=2}^{r} \left( A_{1j}^{\alpha}\right)^{2} \\&& + 2 \sum \limits_{2 \leq i<j\leq r} R^{\left(\ker F_{\ast}\right)^{\perp}} \left(h_{i}, h_{j}, h_{j}, h_{i} \right)- 2 \sum \limits_{2 \leq i < j \leq r} R^{M_{1}} \left(h_{i}, h_{j}, h_{j}, h_{i}\right). 
		\end{eqnarray*}
		Therefore, it follows that
		\begin{equation*}
			2 ~{\rm Ric}_{{\cal H}}^{\left(\ker F_{\ast }\right)^{\perp}} (h_{1}) = 2~ {\rm Ric}_{{\cal H}}^{M_{1}} (h_{1}) - 6 \sum \limits_{\alpha =1}^{n} \sum \limits_{j=2}^{m} \left(A_{1j}^{\alpha }\right)^{2}.
		\end{equation*}
		This finally gives $(\ref{eq-GFCRH-(1)})$ and its equality case.\\
		
		Furthermore, by \cite[Equation (4.11)]{Aytimur_Ozgur}, we have the scalar curvature $\tau^{M_{1}}$ of $M_{1}$ 
		\begin{align}\label{eq-(5.1)}
			\tau^{M_{1}} =& \sum \limits_{1 \leq i<j \leq n} R^{M_{1}} \left(\vee_{i},\vee_{j},\vee_{j},\vee_{i}\right) + \sum \limits_{1 \leq i<j\leq r} R^{M_{1}} \left(h_{i},h_{j},h_{j},h_{i}\right) \nonumber \\&+ \sum \limits_{i=1}^{r} \sum \limits_{j=1}^{n} R^{M_{1}} \left(h_{i},\vee_{j},\vee_{j},h_{i}\right).
		\end{align}
		On the other hand, using (\ref{eq-(2.3)}), (\ref{eq-(2.4)}) and (\ref{eq-(2.5)}), one can get 
		\begin{eqnarray*}
			2\tau^{M_{1}} &=&2\tau_{{\cal H}}^{\left( \ker F_{\ast }\right)^{\perp}}+2\tau_{{\cal V}}^{\ker F_{\ast }}+n^{2}\left\Vert H\right\Vert^{2}\\&&+3 \sum \limits_{i,j=1}^{r}g_{1}\left(A_{h_{i}}h_{j},A_{h_{i}}h_{j}\right)- \sum \limits_{t=1}^{r} \sum \limits_{i,j=1}^{n}\left( T_{ij}^{t}\right)^{2} \\&&-2 \sum \limits_{j=1}^{n} \sum \limits_{i=1}^{r}g_{1}\left( \left( \nabla_{h_{i}}^{1}T\right) \left( \vee_{j},\vee_{j}\right), h_{i}\right) \\&&+2 \sum \limits_{j=1}^{n} \sum \limits_{i=1}^{r}\left\{ g_{1}\left(T_{\vee_{j}}h_{i},T_{\vee_{j}}h_{i}\right) -g_{1}\left(A_{h_{i}}\vee_{j},A_{h_{i}}\vee_{j}\right) \right\}.
		\end{eqnarray*}
		Utilizing (\ref{eq-(3.3)}), (\ref{TVert}), (\ref{AHor}), and (\ref{eq-(3.5)}), we rewrite the aforementioned equation as
		\begin{align}\label{eq-(5.2)}
			\tau^{M_{1}}& = \sum \limits_{1\leq i<j\leq n}R^{\ker F_{\ast }}\left(\vee_{i},\vee_{j},\vee_{j},\vee_{i}\right) + \sum \limits_{1\leq i<j\leq r}R^{\left( \ker F_{\ast}\right)^{\perp }}\left( h_{i},h_{j},h_{j},h_{i}\right) +\frac{n^{2}}{4} \left\Vert H\right\Vert^{2} \nonumber \\& \quad -\frac{1}{4} \sum \limits_{t=1}^{r}\left[ T_{11}^{t}-T_{22}^{t}-\cdots-T_{nn}^{t}\right]^2 - \sum \limits_{t=1}^{r} \sum \limits_{j=2}^{n}\left( T_{1j}^{t}\right)^{2}+ \sum \limits_{t=1}^{r} \sum \limits_{2\leq i<j\leq n}\left\{T_{ii}^{t}T_{jj}^{t}-\left(T_{ij}^{t}\right)^{2}\right\} \nonumber \\& \quad +3 \sum \limits_{\alpha =1}^{n} \sum \limits_{j=2}^{r}\left( A_{1j}^{\alpha }\right)^{2}+3 \sum \limits_{\alpha =1}^{n} \sum \limits_{2\leq i<j\leq r}\left( A_{ij}^{\alpha}\right)^{2}-\delta \left( N\right) +\left\Vert T^{{\cal V}}\right\Vert^{2}-\left\Vert A^{{\cal H}}\right\Vert^{2}.
		\end{align}
		Using (\ref{eq-(3.6a)}) and (\ref{eq-(4.3)}) in (\ref{eq-(5.2)}), we obtain 
		\begin{align}\label{tau_M1}
			\tau^{M_{1}} =& \sum \limits_{1\leq i<j\leq n}R^{\ker F_{\ast }}\left(\vee_{i},\vee_{j},\vee_{j},\vee_{i}\right) + \sum \limits_{1\leq i<j\leq r}R^{\left( \ker F_{\ast}\right)^{\perp }}\left( h_{i},h_{j},h_{j},h_{i}\right) \nonumber \\&-\frac{1}{4} \sum \limits_{t=1}^{r}\left[ T_{11}^{t}-T_{22}^{t}-\cdots -T_{nn}^{t}	\right]^2 - \sum \limits_{t=1}^{r} \sum \limits_{j=2}^{n}\left( T_{1j}^{t}\right)^{2} +3 \sum \limits_{\alpha =1}^{n} \sum \limits_{j=2}^{r}\left(A_{1j}^{\alpha }\right)^{2} \nonumber \\&+ \sum \limits_{2\leq i<j\leq n}R^{M_{1}}\left( \vee_{i},\vee_{j},\vee_{j},\vee_{i}\right)- \sum \limits_{2\leq i<j\leq n}R^{\ker F_{\ast }}\left(\vee_{i},\vee_{j},\vee_{j},\vee_{i}\right) \nonumber \\&+ \sum \limits_{2\leq i<j\leq r}R^{M_{1}}\left( h_{i},h_{j},h_{j},h_{i}\right)- \sum \limits_{2\leq i<j\leq r}R^{\left( \ker F_{\ast }\right)^{\perp }}\left(h_{i},h_{j},h_{j},h_{i}\right) \nonumber \\&-\delta \left( N\right) +\left\Vert T^{{\cal V}}\right\Vert^{2}-\left\Vert A^{{\cal H}}\right\Vert^{2}+\frac{n^{2}}{4}\left\Vert H\right\Vert^{2}.
		\end{align}
		Then by (\ref{eq-(5.1)}) and (\ref{tau_M1}), we have 
		\begin{align*}
			{\rm Ric}_{{\cal V}}^{M_{1}}\left( \vee_{1}\right) &+{\rm Ric}_{{\cal H}}^{M_{1}}\left( h_{1}\right) + \sum \limits_{i=1}^{r} \sum \limits_{j=1}^{n}R^{M_{1}}\left(h_{i},\vee_{j},\vee_{j},h_{i}\right) \\&={\rm Ric}^{\ker F_{\ast }}_{\cal V}\left( \vee_{1}\right) +{\rm Ric}^{\left( \ker	F_{\ast }\right)^{\perp }}_{\cal H}\left( h_{1}\right) -\frac{1}{4} \sum \limits_{t=1}^{r}\left[ T_{11}^{t}-T_{22}^{t}-\cdots -T_{nn}^{t}\right]^2 \\&- \sum \limits_{t=1}^{r} \sum \limits_{j=2}^{n}\left( T_{1j}^{t}\right)^{2}+3 \sum \limits_{\alpha=1}^{n} \sum \limits_{j=2}^{r}\left( A_{1j}^{\alpha}\right)^{2}-\delta \left(N\right) +\left\Vert T^{{\cal V}}\right\Vert^{2}-\left\Vert A^{{\cal H}}\right\Vert^{2}+\frac{n^{2}}{4}\left\Vert H\right\Vert^{2}.
		\end{align*}
		This implies $(\ref{eq-GFCRVH-(1)})$ and its equality case.
	\end{proof}
	
	\subsection{Example}
	Here, we provide an example satisfying the equality of Chen-Ricci inequalities derived in Theorem \ref{Theorem-GCRV}, as follows.
	
	\begin{example}
		Considering $\left(M_1 = (- \frac{\pi}{2}, \frac{\pi}{2}) \times_{\cos t} \mathbb{C}, ~g_1 = dt^2 + \cos^2 t \sum \limits_{i=1}^{3} ((dx_i)^2 + (dy_i)^2)\right)$ and $\left( M_2 =\mathbb{R}^4, g_2 = \cos^2 t \sum \limits_{j=1}^{4} (dz_j)^2\right)$, we define $F: (M_1, g_1) \to (M_2, g_2)$ such that $$F(t, x_1, y_1, x_2, y_2, x_3, y_3) = (x_2, y_2, x_3, y_3).$$ Then we get
		\begin{equation*}
			(\ker F_\ast) = \operatorname{span}\left\{\vee_1 = \frac{\partial}{\partial t}, \vee_2 = \frac{1}{\cos t} \frac{\partial}{\partial x_1}, \vee_3 = \frac{1}{\cos t} \frac{\partial}{\partial y_1}\right\},
		\end{equation*}
		\begin{equation*}
			(\ker F_\ast)^\perp = \operatorname{span}\left\{h_1 = \frac{1}{\cos t} \frac{\partial}{\partial x_2}, h_2 = \frac{1}{\cos t} \frac{\partial}{\partial y_2}, h_3 = \frac{1}{\cos t} \frac{\partial}{\partial x_3}, h_4 = \frac{1}{\cos t} \frac{\partial}{\partial y_3}\right\},
		\end{equation*}
		and
		\begin{equation*}
			({\rm range}~ F_\ast) = \operatorname{span}\left\{F_\ast (h_j) = \frac{1}{\cos t} \frac{\partial}{\partial z_j} \right\}_{j=1}^{4},
		\end{equation*}
		where $\left\{\frac{\partial}{\partial t}, \frac{1}{\cos t} \frac{\partial}{\partial x_1}, \frac{1}{\cos t} \frac{\partial}{\partial y_1}, \frac{1}{\cos t} \frac{\partial}{\partial x_2}, \frac{1}{\cos t} \frac{\partial}{\partial y_2}, \frac{1}{\cos t} \frac{\partial}{\partial x_3}, \frac{1}{\cos t} \frac{\partial}{\partial y_3} \right\}$ and $\left\{\frac{\partial}{\partial z_j} \right\}_{j=1}^{4}$ are bases of $T_p M_1$ and $T_{F(p)}M_2$, respectively. For all $i, j \in \{1,2,3,4\}$, we verify that $$g_1 (h_i, h_j) = g_2 (F_\ast h_i, F_\ast h_j).$$ Therefore, $F$ is a Riemannian submersion between Riemannian manifolds.\\
		
		We calculate the non-vanishing Christoffel symbols for the metric $g_1$, which are $\Gamma^t_{i, i} = \sin t \cos t$ and $\Gamma^i_{t, i} = -\tan t$ for all $i \in \{x_j, y_j\}_{j=1}^{3}$. Using these values, we calculate the following:
		\begin{equation*}
			\nabla^{M_1}_{\vee_1} \vee_i = \nabla^{M_1}_{\vee_2} \vee_3 = \nabla^{M_1}_{\vee_3} \vee_2 = 0 ~ \text{for all $1 \leq i \leq 3$},
		\end{equation*}
		\begin{equation*}
			\nabla^{M_1}_{\vee_2} \vee_1 = -\tan t \sec t\frac{\partial}{\partial x_1}, \nabla^{M_1}_{\vee_3} \vee_1 = -\tan t \sec t \frac{\partial}{\partial y_1}, \nabla^{M_1}_{\vee_2} \vee_2 = \nabla^{M_1}_{\vee_3} \vee_3 = \tan t \frac{\partial}{\partial t},
		\end{equation*}
		\begin{equation*}
			\nabla^{M_1}_{h_i} h_j = 0 ~ \text{for all $1  \leq i \neq j \leq 4$ and ~} \nabla^{M_1}_{h_i} h_i = \tan t \frac{\partial}{\partial t} ~ \text{for all $1  \leq i \leq 4$}.
		\end{equation*}
		Thus, $\mathcal{H} \nabla^{M_1}_{\vee_i} \vee_j = T_{\vee_i} \vee_j = 0$ for all $1 \leq i, j \leq 3$. Hence, for all $h_t \in \{h_1, h_2, h_3, h_4\}$, we have $$T_{ij}^{t}=g_1(T_{\vee_i} \vee_j, h_t) = 0, \quad 1 \leq i, j \leq 3 ~ \text{and ~} 1 \leq t \leq 4.$$
		In addition, we see that $A_{h_1} h_j =0$ for all $j \in \{2, 3, 4\}$. Hence $A_{1j}^{\alpha} = g_1 (A_{h_1} h_j, \vee_\alpha) = 0$ for all $j \in \{2, 3, 4\}$ and $\alpha \in \{1, 2, 3\}$. Consequently, for $F$ the inequalities mentioned in $(\ref{eq-GFCRV-(1)})$, $(\ref{eq-GFCRH-(1)})$ and $(\ref{eq-GFCRVH-(1)})$ attain equality.
	\end{example}
	
	\subsection{Applications}\label{appl_for_part_I}
	
	In this subsection, we obtain the Chen-Ricci inequalities for Riemannian submersions from generalized complex and generalized Sasakian space forms. Toward this, we use the general forms established in Theorem \ref{Theorem-GCRV}. Also, for any ${\cal Y}\in \Gamma (TM_1)$, we decompose
	\begin{equation*}
		J {\cal Y}  {\textrm ~or~ } \phi {\cal Y}= {Q} {\cal Y}+{P}{\cal Y}, 
	\end{equation*}
	where ${Q} {\cal Y}\in \Gamma ({\rm ker~}F_\ast)$, ${P}{\cal Y}\in \Gamma ({\rm ker~} F_\ast)^{\perp}$ with $ \left \Vert P \right \Vert^{2} = \sum \limits_{i=1}^{r} \sum \limits_{j=1}^{n} \left( g_1\left( P\vee_{j}, h_{i}\right) \right)^{2}$, $\left \Vert P h_{1}\right\Vert^{2} = \sum \limits_{j=2}^{r} \left( g_1\left( P h_1, h_{j}\right) \right)^{2}$, and $\left \Vert Q \vee_{1}\right\Vert^{2} = \sum \limits_{j=2}^{n}\left( g_1 \left(Q\vee_{1},\vee_{j}\right) \right)^{2}$.
	
	\begin{theorem}\label{appl_for_CRI_RS}
		Suppose $F$ satisfies the hypotheses of Theorem $\ref{Theorem-GCRV}$. Then the following hold.
		\begin{enumerate}[$(i)$]
			\item If $(M_1(c_1, c_2), g_1, J)$ is a generalized complex space form, then
			\begin{equation}\label{eq-GCCRV-(1)}
				{\rm Ric}_{{\cal V}}^{\ker F_{\ast }}(\vee_{1})\geq c_{1}\left( n-1\right) + 3 c_{2} \left\Vert Q\vee_{1}\right\Vert^{2} - \frac{1}{4}n^{2}\left\Vert H\right\Vert^{2},
			\end{equation}
			\begin{equation}\label{eq-GCCRH-(1)}
				{\rm Ric}_{{\cal H}}^{\left( \ker F_{\ast }\right)^{\perp }}(h_{1})\leq c_{1}\left( r-1\right) +3c_{2}\left\Vert Ph_{1}\right\Vert^{2},
			\end{equation}
			and
			\begin{align}\label{eq-GCCRVH-(1)}
				{\rm Ric}_{{\cal V}}^{\ker F_\ast}\left( \vee_{1}\right) &+{\rm Ric}_{\cal H}^{(\ker F_\ast)^\perp}\left(h_{1}\right)\nonumber\\& \geq c_{1}\left\{ nr+n+r-2\right\} +3c_{2}\left\{ \left\Vert P\right\Vert^{2}+\left\Vert Ph_{1}\right\Vert^{2}+\left\Vert Q\vee_{1}\right\Vert^{2}\right\} \nonumber \\&+\delta (N) - \left\Vert T^{{\cal V}}\right\Vert^{2} + \left\Vert A^{{\cal H}}\right\Vert^{2} - \frac{n^{2}}{4}\left\Vert H\right\Vert^{2} - 3 \sum \limits_{\alpha=1}^{n} \sum \limits_{j=2}^{r}\left( A_{1j}^{\alpha }\right)^{2}.
			\end{align}
			
			\item If $\left({M_{1}} \left( c_{1},c_{2},c_{3}\right), \phi, \xi, \eta, g_{1}\right)$ is a generalized Sasakian space form, then
			\begin{align}\label{eq-VGSCRV-(1)}
				{\rm Ric}_{{\cal V}}^{\ker F_{\ast }}(\vee_{1}) - c_{1}\left( n-1\right) &+ \frac{1}{4} n^{2} \left\Vert H\right \Vert^{2} - 3c_{2}\left\Vert Q\vee_{1}\right\Vert^{2} \nonumber\\& \geq \left\{\begin{array}{l}
					- c_{3}\left(1+\left( n-2\right) \left( \eta \left( \vee_{1}\right) \right)^{2}\right), ~{\rm if}~ \xi ~\text{is vertical}; \\ 0, ~{\rm if}~ \xi ~\text{is horizontal},
				\end{array} \right.
			\end{align}
			\begin{align}\label{eq-VGSCRH-(1)}
				{\rm Ric}_{{\cal H}}^{(\ker F_{\ast})^\perp}(h_{1}) &- c_{1} \left( r-1\right) - 3c_{2} \left \Vert P h_{1}\right\Vert^{2} \nonumber\\& \leq \left\{\begin{array}{l}
					0, ~{\rm if}~ \xi ~\text{is vertical}; \\ -c_{3}\left(1+\left( r-2\right) \left( \eta \left( h_{1}\right) \right)^{2}\right), ~{\rm if}~ \xi ~\text{is horizontal},
				\end{array} \right.
			\end{align}
			and
			\begin{align}\label{eq-VGSCRVH-(1)}
				& {\rm Ric}_{{\cal V}}^{\ker F_\ast}\left( \vee_{1}\right) +{\rm Ric}_{\cal H}^{(\ker F_\ast)^\perp}\left(h_{1}\right) \nonumber \\&+\frac{n^{2}}{4}\left\Vert H\right\Vert^{2}+3 \sum \limits_{\alpha=1}^{n} \sum \limits_{j=2}^{r}\left( A_{1j}^{\alpha }\right)^{2} -\delta \left( N\right) +\left\Vert T^{{\cal V}}\right\Vert^{2}-\left\Vert A^{{\cal H}}\right\Vert^{2} \nonumber \\&- c_{1}\left\{ nr+n+r-2\right\} - 3c_{2}\left\{ \left\Vert P\right\Vert^{2}+\left\Vert Ph_{1}\right\Vert^{2}+\left\Vert Q\vee_{1}\right\Vert^{2}\right\}\nonumber \\&\geq \left\{\begin{array}{l}
					- c_{3}\left\{ r+1+\left( n-2\right) \left( \eta \left( \vee_{1}\right)\right)^{2}\right\}, ~{\rm if}~ \xi ~\text{is vertical}; \\ -c_{3}\left\{ n+1+\left( r-2\right) \left( \eta \left( h_{1}\right)\right)^{2}\right\}, ~{\rm if}~ \xi ~\text{is horizontal}.
				\end{array} \right.
			\end{align}
		\end{enumerate} 
		
		In addition, the equality cases are similar to those of Theorem $\ref{Theorem-GCRV}$.
	\end{theorem}

	\begin{proof}[Proof of $(i)$]
		Using (\ref{curvature_for_gcsf}), we derive
		\begin{equation}\label{eq-LGCRV-(1)}
			{\rm Ric}_{{\cal V}}^{M_{1}} \left( \vee_{1}\right) = c_{1} \left( n-1\right) + 3 c_{2} \left \Vert Q\vee_{1} \right \Vert^{2},
		\end{equation}
		\begin{equation}\label{eq-LGCRH-(1)}
			{\rm Ric}_{{\cal H}}^{M_{1}}\left( h_{1}\right) =c_{1}\left( r-1\right) + 3 c_{2}\left\Vert Ph_{1}\right\Vert^{2},
		\end{equation}
		and
		\begin{equation}\label{eq-LGCRVH-(1)}
			\sum \limits_{i=1}^{r} \sum \limits_{j=1}^{n}R^{M_{1}}\left(\vee_{j},h_{i},h_{i},\vee_{j}\right) = nrc_{1} + 3 c_{2} \Vert P \Vert^{2}.
		\end{equation}
		Adding equations (\ref{eq-LGCRV-(1)}), (\ref{eq-LGCRH-(1)}), and (\ref{eq-LGCRVH-(1)}), we obtain
		\begin{align}\label{eq-LGCRVH-(2)}
			{\rm Ric}_{{\cal V}}^{M_{1}}\left( \vee_{1}\right) &+{\rm Ric}_{{\cal H}}^{{M}_1}\left( h_{1}\right) + \sum \limits_{i=1}^{r} \sum \limits_{j=1}^{n}R^{M_{1}}\left(\vee_{j},h_{i},h_{i},\vee_{j}\right) \nonumber \\& = c_{1} \left\{ nr+n+r-2\right\} +3c_{2}\left\{\left\Vert P\right\Vert^{2}+\left\Vert Ph_{1}\right\Vert^{2}+\left\Vert Q\vee_{1}\right\Vert^{2}\right\}.
		\end{align}
		Then by using (\ref{eq-LGCRV-(1)}) and (\ref{eq-GFCRV-(1)}), (\ref{eq-LGCRH-(1)}) and (\ref{eq-GFCRH-(1)}), (\ref{eq-LGCRVH-(2)}) and (\ref{eq-GFCRVH-(1)}), we obtain the desired inequalities (\ref{eq-GCCRV-(1)}), (\ref{eq-GCCRH-(1)}) and (\ref{eq-GCCRVH-(1)}) respectively.
	\end{proof}
	
	\begin{proof}[Proof of $(ii)$]
		Using (\ref{curvature_for_gssf}), we derive
		\begin{align}\label{(1)}
			&{\rm Ric}_{{\cal V}}^{M_{1}}\left( \vee_{1}\right) = \left\{\begin{array}{l}
				c_{1}\left( n-1\right) + 3 c_{2} \left\Vert Q\vee_{1} \right \Vert^{2} - c_{3} \left( 1+\left( n-2\right) \left( \eta \left( \vee_{1}\right) \right)^{2}\right), ~{\rm if}~ \xi ~\text{is vertical}; \\ c_{1}\left( n-1\right) + 3 c_{2} \left \Vert Q\vee_{1}\right\Vert^{2}, ~{\rm if}~ \xi ~\text{is horizontal},
			\end{array} \right.
		\end{align}
		\begin{align}\label{(2)}
			&{\rm Ric}_{{\cal H}}^{M_{1}}\left( h_{1}\right) = \left\{\begin{array}{l}
				c_{1}\left( r-1\right) + c_{2}\left\Vert Ph_{1}\right\Vert^{2}, ~{\rm if}~ \xi ~\text{is vertical}; \\ c_{1}\left( r-1\right) + 3 c_{2}\left\Vert Ph_{1}\right\Vert^{2}-c_{3}\left( 1+\left( r-2\right) \left( \eta \left( h_{1}\right) \right)^{2}\right), ~{\rm if}~ \xi ~\text{is horizontal},
			\end{array} \right.
		\end{align}
		and
		\begin{align}\label{(3)}
			&\sum \limits_{i=1}^{r} \sum \limits_{j=1}^{n}R^{M_{1}}\left( \vee_{j},h_{i},h_{i},\vee_{j}\right) = \left\{\begin{array}{l}
				nrc_{1} +3c_{2} \Vert P \Vert^{2}-c_{3}r, ~{\rm if}~ \xi ~\text{is vertical}; \\ nrc_{1} +3c_{2} \Vert P\Vert^{2}-c_{3}n, ~{\rm if}~ \xi ~\text{is horizontal}.
			\end{array} \right.
		\end{align}
		Adding equations (\ref{(1)}), (\ref{(2)}), and (\ref{(3)}), we obtain
		\begin{align}\label{(4)}
			&{\rm Ric}_{{\cal V}}^{M_{1}}\left( \vee_{1}\right) +{\rm Ric}_{{\cal H}}^{{M}_1}\left( h_{1}\right) \nonumber \\&+ \sum \limits_{i=1}^{r} \sum \limits_{j=1}^{n}R^{M_{1}}\left(\vee_{j},h_{i},h_{i},\vee_{j}\right)-c_{1}\left\{ nr+n+r-2\right\} \nonumber \\&-3c_{2}\left\{ \left\Vert P\right\Vert^{2}+\left\Vert Ph_{1}\right\Vert^{2}+\left\Vert Q\vee_{1}\right\Vert^{2}\right\} \nonumber \\& = \left\{\begin{array}{l}
				-c_3 \left\{ r+1+\left( n-2\right) \left( \eta \left(\vee_{1}\right) \right)^{2}\right\}, ~{\rm if}~ \xi ~\text{is vertical}; \\ -c_{3}\left\{ n+1+\left( r-2\right) \left( \eta \left( h_{1}\right) \right)^{2}\right\}, ~{\rm if}~ \xi ~\text{is horizontal}.
			\end{array} \right.
		\end{align}
		Then by using (\ref{(1)}) and (\ref{eq-GFCRV-(1)}), (\ref{(2)}) and (\ref{eq-GFCRH-(1)}), (\ref{(4)}) and (\ref{eq-GFCRVH-(1)}), we obtain the desired inequalities (\ref{eq-VGSCRV-(1)}), (\ref{eq-VGSCRH-(1)}) and (\ref{eq-VGSCRVH-(1)}) respectively.
	\end{proof}
	
	\section{General Chen-Ricci inequalities for Riemannian maps with applications}\label{part_II}
	
	\subsection{Background}\label{sec_2}
	
	In this subsection, we recall some information required for the next subsections.\\
	
	Let $F:(M_{1}^{m_{1}},g_{1}) \to (M_{2}^{m_{2}},g_{2})$ be a smooth map between Riemannian manifolds with rank $r< \min \{m_1, m_2\}$, and let $F_{\ast p} : T_{p} M_{1}\to T_{F(p)} M_{2}$ be its derivative map at $p$. Denoting the kernel (resp. range) space of $F_{\ast}$ by $(\ker F_{\ast p})$ (resp. $({\rm range~} F_{\ast p})$) and its orthogonal complementary space by $(\ker F_{\ast p})^{\perp}$ (resp. $({\rm range~} F_{\ast p})^\perp$), we decompose 
	\[
	T_{p}M_{1}=(\ker F_{\ast p})\oplus (\ker F_{\ast p})^{\perp } ~{\rm and~} T_{F(p)}M_{2}=({\rm range~}F_{\ast p})\oplus ({\rm range~}F_{\ast p})^{\perp}.
	\]
	The map $F$ is called a {\it Riemannian map}, if for all $X,Y \in \Gamma (\ker F_{\ast })^{\perp }$, we have \cite{Fischer_1992}
	\begin{equation}\label{eq-(1)}
		g_{1}(X,Y)=g_{2}(F_{\ast }X,F_{\ast }Y) .
	\end{equation}
	
	\subsubsection*{Second fundamental form}
	\noindent The bundle {\rm Hom}$\left( T{M_{1}},F^{-1}T{M_{2}}\right)$ admits an induced connection $\nabla$ from the Levi-Civita connection $\nabla^{1}$ on ${M_{1}}$. The symmetric {\it second fundamental form} of $F$ is then given by \cite{Nore_1986}
	\begin{equation}\label{2ndff}
		\left( \nabla F_{\ast }\right) \left( {\cal Z}_{1},{\cal Z}_{2}\right) =\nabla_{F_\ast {\cal Z}_{1}}^{2}F_{\ast }({\cal Z}_{2})-F_{\ast }\left( \nabla_{{\cal Z}_{1}}^{{{1}}}{\cal Z}_{2}\right)
	\end{equation}
	for ${\cal Z}_{1},{\cal Z}_{2}\in \Gamma (T{M_{1}})$. In addition, by \cite{Sahin_2010} $\left( \nabla F_{\ast }\right) \left( X,Y\right) $ is completely contained in $\left( {\rm range~}F_{\ast }\right)^{\perp }$ for all $X,Y\in \Gamma \left( \ker F_{\ast }\right)^{\perp }$. 
	
	\subsubsection*{Gauss equation}
	The \textit{Gauss equation} for $F$ is defined as \cite[p. 189]{Sahin_book}
	\begin{eqnarray}\label{Gauss_Eqn}
		g_{2}\left( R^{{M_{2}}}\left( F_{\ast }Z_{1},F_{\ast }Z_{2}\right) F_{\ast}Z_{3},F_{\ast }Z_{4}\right) &=&g_{1}\left( R^{{M_{1}}}\left(Z_{1},Z_{2}\right) Z_{3},Z_{4}\right)  \nonumber \\&&+g_{2}\left( \left( \nabla F_{\ast }\right) \left( Z_{1},Z_{3}\right),\left( \nabla F_{\ast }\right) \left(Z_{2},Z_{4}\right) \right)  \nonumber\\&&-g_{2}\left( \left( \nabla F_{\ast }\right) \left( Z_{1},Z_{4}\right),\left( \nabla F_{\ast }\right) \left( Z_{2},Z_{3}\right) \right) ,
	\end{eqnarray}
	where $Z_{i}\in \Gamma \left( \ker F_{\ast }\right)^{\perp }$. Here, $R^{M_{1}}$ and $R^{M_{2}}$ denote the curvature tensors of ${M_{1}}$ and ${M_{2}}$, respectively.
	
	\subsubsection*{Some notations}
	
	At a point $p\in {M_{1}}$, suppose that $\{h_{i}\}_{i=1}^{r}$ and $\{h_{i}\}_{i=r+1}^{m_{1}}$ are orthonormal bases of the horizontal space (that is, $(\ker F_\ast)^\perp$) and the vertical space (that is, $(\ker F_\ast)$), respectively. Then the scalar curvatures $\tau^{{\cal H}}$ and $\tau^{{\cal R}}$ on the horizontal and range spaces are given, respectively, by 
	\[
	\tau^{{\cal H}}=\sum\limits_{1\leq i<j\leq r}g_{1}\left( R^{{M_{1}}}(h_{i},h_{j})h_{j},h_{i}\right),\quad \tau^{{\cal R}}=\sum\limits_{1\leq i<j\leq r}g_{2}\left( R^{{M_{2}}}(F_{\ast }h_{i},F_{\ast }h_{j})F_{\ast }h_{j},F_{\ast }h_{i}\right). 
	\]
	In addition, for a fixed $j$, the Ricci curvatures of $h_{j}$ and $F_{\ast }h_{j}$, denoted by ${\rm Ric}(h_{j})$ and ${\rm Ric}(F_{\ast }h_{j})$, respectively, are defined by
	\begin{equation}\label{Ricci_Curvatures}
		{\rm Ric}^{\cal H}(h_{j}) =\sum_{i= 1}^{r}g_{1}\left( R^{{M_{1}}}(h_{i},h_{j})h_{j},h_{i}\right),  {\rm Ric}^{\cal R}(F_{\ast }h_{j}) =\sum_{i =1}^{r}g_{2}\left( R^{{M_{2}}}(F_{\ast }h_{i},F_{\ast }h_{j})F_{\ast }h_{j},F_{\ast }h_{i}\right).
	\end{equation}
	Supposing $\{V_{r+1}, \dots, V_{m_2}\}$ an orthonormal basis of $\left( {\rm range~}F_{\ast}\right)^{\perp }$ we fix, 
	\begin{eqnarray*}
		B_{ij}^{{\cal H}^{\alpha }} &=&g_{2}\left( (\nabla F_{\ast})(h_{i}, h_{j}), V_{\alpha }\right), \quad 1 \leq i, j \leq r, \quad r+1 \leq \alpha \leq m_2, \nonumber \\ {\rm trace\, }B^{{\cal H}} &=&\sum_{i=1}^{r}(\nabla F_{\ast}) \left(h_{i}, h_{i}\right), \quad {\rm and} \nonumber \\ \left\Vert {\rm trace\, }B^{{\cal H}}\right\Vert^{2} &=&g_{2}\left( {\rm trace\, }B^{{\cal H}}, {\rm trace\, }B^{{\cal H}}\right).
	\end{eqnarray*}
	Then we get \cite{LLSV_2022},
	\begin{eqnarray}\label{eq-(5)}
		\left\Vert B^{{\cal H}}\right\Vert^{2} &=&\sum_{i, j=1}^{r}g_{2}\left((\nabla F_{\ast})(h_{i}, h_{j}), (\nabla F_{\ast})(h_{i}, h_{j})\right) \nonumber \\ &=&\frac{1}{2}\left\Vert {\rm trace~}	B^{{\cal H}}\right\Vert^{2}+\frac{1}{2}\sum_{\alpha =r+1}^{m_2}\left( B_{11}^{{\cal H}^{\alpha }}-B_{22}^{{\cal H}^{\alpha }}-\cdots -B_{rr}^{{\cal H}^{\alpha }}\right)^{2}  \nonumber \\&&+2\sum_{\alpha =r+1}^{m_2}\sum_{i=2}^{r}\left( B_{1i}^{{\cal H}^{\alpha}}\right)^{2}-2\sum_{\alpha =r+1}^{m_2}\sum_{2\leq i<j\leq r}\left\{ B_{ii}^{{\cal H}^{\alpha }}B_{jj}^{{\cal H}^{\alpha }}-\left( B_{ij}^{{\cal H}^{\alpha }}\right)^{2}\right\}.  
	\end{eqnarray}

	\subsection{General Chen-Ricci inequality}\label{sec_3}
	
	Here, we present general form of the Chen-Ricci inequality for Riemannian maps between Riemannian manifolds. Toward this, we have the following result.
	
	\begin{theorem}\label{main_thm_CRI}
		Let $F:\left(M_{1}^{m_1},g_{1}\right) \to \left( M_{2}^{m_2},g_{2}\right) $ be a Riemannian map between Riemannian manifolds with rank $r<m_2$. Then for any unit vector $X \in \Gamma (\ker F_\ast)^\perp$, we have
		\begin{equation}\label{eq-(3)}
			{\rm Ric}^{{\cal H}}(X) \leq {\rm Ric}^{{\cal R}}(F_{\ast }X) + \frac{1}{4} \left\Vert {\rm trace~}B^{{\cal H}}\right\Vert^{2}.
		\end{equation}
		\begin{enumerate}[$(i)$]
			\item The equality case of $(\ref{eq-(3)})$ holds for all $X\in \Gamma (\ker F_\ast)^\perp$ if and only if
			\begin{eqnarray*}
				\left( \nabla F_{\ast }\right) \left( X,Y\right) &=&0,\quad {\rm for\ all\ } Y\in \Gamma (\ker F_\ast)^\perp \ {\rm orthogonal\ to\ }X,  \nonumber \\ \left( \nabla F_{\ast }\right) \left( X,X\right) &=&\frac{{\rm trace~}B^{{\cal H}}}{2}.
			\end{eqnarray*}
			
			\item The equality case of $(\ref{eq-(3)})$ holds identically for all $X\in \Gamma (\ker F_\ast)^\perp$ if and only if either $B^{{\cal H}}=0$ or $r=2$, and $B_{11}^{{\cal H}^{\alpha }}=B_{22}^{{\cal H}^{\alpha }}$, $\alpha \in \{3, \ldots, m_2\}$. Equivalently, the equality case of $(\ref{eq-(3)})$ holds identically for all $X\in \Gamma (\ker F_\ast)^\perp$ if and only if $F$ is
			\begin{enumerate}[$(a)$]
				\item totally geodesic when $r>2$;
				
				\item umbilical when $r=2$.
			\end{enumerate}
		\end{enumerate}
	\end{theorem}
	
	\begin{proof}
		Substituting $Z_{1}=Z_{4}=h_{i}$ and $Z_{2}=Z_{3}=h_{j} $ in (\ref{Gauss_Eqn}), we obtain
		\begin{equation}\label{eq-(4)}
			2\tau^{{\cal H}}=2\tau^{{\cal R}}+\left\Vert {\rm trace~}B^{{\cal H}}\right\Vert^{2}-\left\Vert B^{{\cal H}}\right\Vert^{2}. 
		\end{equation}
		Employing (\ref{eq-(5)}) in (\ref{eq-(4)}), we get
		\begin{eqnarray}\label{eq-(6)}
			\tau^{{\cal H}} &=&\tau^{{\cal R}}+\frac{1}{4}\left\Vert {\rm trace~}B^{				{\cal H}}\right\Vert^{2}-\frac{1}{4}\sum_{\alpha =r+1}^{m_2}\left( B_{11}^{				{\cal H}^{\alpha }}-B_{22}^{{\cal H}^{\alpha }}-\cdots -B_{rr}^{{\cal H}^{\alpha }}\right)^{2}  \nonumber \\	&&-\sum_{\alpha =r+1}^{m_2}\sum_{i=2}^{r}\left( B_{1i}^{{\cal H}^{\alpha}}\right)^{2}+\sum_{\alpha =r+1}^{m_2}\sum_{2\leq i<j\leq r}\left\{ B_{ii}^{{\cal H}^{\alpha }}B_{jj}^{{\cal H}^{\alpha }}-\left( B_{ij}^{{\cal H}^{\alpha }}\right)^{2}\right\}.
		\end{eqnarray}
		By (\ref{Gauss_Eqn}), we have
		\begin{eqnarray}\label{eq-(7)}
			&&\sum_{\alpha =r+1}^{m_2}\sum_{2\leq i<j\leq r}\left\{ B_{ii}^{{\cal H}^{\alpha }}B_{jj}^{{\cal H}^{\alpha }}-\left( B_{ij}^{{\cal H}^{\alpha
			}}\right)^{2}\right\}  \nonumber \\
			&&=\sum_{2\leq i<j\leq r}g_{1}\left( R^{M_{1}}\left( h_{i},h_{j}\right)
			h_{j},h_{i}\right) -\sum_{2\leq i<j\leq r}g_{2}\left( R^{M_{2}}\left(
			F_{\ast }h_{i},F_{\ast }h_{j}\right) F_{\ast }h_{j},F_{\ast }h_{i}\right).
		\end{eqnarray}
		Then substituting (\ref{eq-(7)}) in (\ref{eq-(6)}), we get
		\begin{eqnarray*}
			\tau^{{\cal H}}-&&\sum_{2\leq i<j\leq r}g_{1}\left( R^{M_{1}}\left(h_{i},h_{j}\right) h_{j},h_{i}\right)\\&&=\tau^{{\cal R}}-\sum_{2\leq i<j\leq r}g_{2}\left( R^{M_{2}}\left(F_{\ast }h_{i},F_{\ast }h_{j}\right) F_{\ast }h_{j},F_{\ast }h_{i}\right) +			\frac{1}{4}\left\Vert {\rm trace~}B^{{\cal H}}\right\Vert^{2} \\			&&-\frac{1}{4}\sum_{\alpha =r+1}^{m_2}\left( B_{11}^{{\cal H}^{\alpha}}-B_{22}^{{\cal H}^{\alpha }}-\cdots -B_{rr}^{{\cal H}^{\alpha }}\right)^{2}-\sum_{\alpha =r+1}^{m_2}\sum_{i=2}^{r}\left( B_{1i}^{{\cal H}^{\alpha}}\right)^{2}.
		\end{eqnarray*}
		Equivalently,
		\begin{eqnarray}\label{eq-(8)}
			{\rm Ric}^{{\cal H}}(h_{1}) &=&{\rm Ric}^{{\cal R}}(F_{\ast }h_{1})+\frac{1}{4}\left\Vert {\rm trace~}B^{{\cal H}}\right\Vert^{2}  \nonumber \\&&-\frac{1}{4}\sum_{\alpha =r+1}^{m_2}\left( B_{11}^{{\cal H}^{\alpha}}-B_{22}^{{\cal H}^{\alpha }}-\cdots -B_{rr}^{{\cal H}^{\alpha }}\right)^{2}-\sum_{\alpha =r+1}^{m_2}\sum_{i=2}^{r}\left( B_{1i}^{{\cal H}^{\alpha}}\right)^{2}.
		\end{eqnarray}
		Then (\ref{eq-(3)}) follows as we can choose $h_{1}=X \in \Gamma \left( \ker F_\ast \right)^\perp$. Also, the equality holds if and only if
		\begin{eqnarray}\label{eq-(9)}
			B_{11}^{{\cal H}^{\alpha }} =B_{22}^{{\cal H}^{\alpha }}+\cdots +B_{rr}^{{\cal H}^{\alpha }}, ~ B_{1i}^{{\cal H}^{\alpha }} =0, \quad i \in \{2,\ldots ,r\},\ \alpha \in \{r+1,\ldots ,m_2\}.
		\end{eqnarray}
		This implies $(i)$.\\
		
		\noindent Further, assume that the equality case of (\ref{eq-(3)}) holds for all $X \in \Gamma (\ker F_\ast)^\perp$. Then by (\ref{eq-(9)}), we have
		\begin{equation}\label{eq-(11)}
			B_{ij}^{{\cal H}^{\alpha }}=0\quad ~{\rm and} \quad 2 B_{ii}^{{\cal H}^{\alpha }}=B_{11}^{{\cal H}^{\alpha }}+B_{22}^{{\cal H}^{\alpha }}+\cdots +B_{rr}^{{\cal H}^{\alpha }},\quad 1 \leq i \neq j \leq r.
		\end{equation}
		This implies
		\begin{eqnarray*}
			\left( r-2\right) \sum \limits_{i=1}^{r} B_{ii}^{{\cal H}^{\alpha }}=0.
		\end{eqnarray*}
		Thus, either $\sum \limits_{i=1}^{r} B_{ii}^{{\cal H}^{\alpha }}=0$ or $r=2$. 
		\begin{itemize}
			\item If $r=2$, then using (\ref{eq-(11)}), we get $B_{11}^{{\cal H}^{\alpha }}=B_{22}^{{\cal H}^{\alpha }}$.
			\item If $\sum \limits_{i=1}^{r} B_{ii}^{{\cal H}^{\alpha }}=0$, then using (\ref{eq-(11)}), we get ${\rm trace~} B^{\cal H} =0$.
		\end{itemize} 
		We also note that the converse case is trivial. Thus, we complete the required proof.
	\end{proof}
	
	\subsection{Improved general Chen-Ricci inequality}\label{sec_4}
	In this subsection, we present general form of the improved Chen-Ricci inequality for Riemannian maps between Riemannian manifolds. Toward this, we have the following result.
	
	\begin{theorem}\label{main_thm_ICRI}
		Let $F:\left(M_{1}^{m_1},g_{1}\right) \to \left( M_{2}^{m_2},g_{2}\right) $ be a Riemannian map between Riemannian manifolds with rank $r < m_2$. Then for any unit vector $X \in \Gamma (\ker F_\ast)^\perp$, we have
		\begin{equation}\label{eq-IC(4.5)}
			{\rm Ric}^{{\cal H}}\left( X\right) \leq {\rm Ric}^{{\cal R}}\left( F_{\ast}X\right) +\frac{(r-1)}{4r}\left\Vert {\rm trace~}B^{{\cal H}}\right\Vert^{2}.
		\end{equation}
		The equality holds for all unit $X \in \Gamma(\ker F_\ast)^\perp$ if and only if either
		\begin{enumerate}[$(i)$]
			\item $B^{{\cal H}}=0$, equivalently, $F$ is totally geodesic, or
			
			\item $r=2$ and $$\left( \nabla F_{\ast }\right) \left( h_{1},h_{1}\right) = 3\mu V_{r+1}, \left( \nabla F_{\ast }\right) \left( h_{2},h_{2}\right) = \mu V_{r+1}, \left( \nabla F_{\ast }\right) \left( h_{1},h_{2}\right) =\mu V_{r+2}$$
			for some function $\mu $ with respect to some orthonormal local frame field $\left\{ h_{1},h_{2}\right\} $.
		\end{enumerate}
	\end{theorem}
	
	\begin{proof}
		For chosen bases, let
		\begin{eqnarray}\label{eq-IC(4.6)}
			B_{ij}^{{\cal H}^{r+k}} =B_{jk}^{{\cal H}^{r+i}}=B_{ki}^{{\cal H}^{r+j}}, B_{ij}^{{\cal H}^{\beta }} = 0, \quad i,j,k\in \left\{ 1,\ldots ,r\right\}, \beta \in \left\{ 2r+1,\ldots,m_2\right\}.
		\end{eqnarray}
		By (\ref{Gauss_Eqn}), we get
		\begin{equation*}
			{\rm Ric}^{{\cal R}}\left( F_{\ast }X\right) ={\rm Ric}^{{\cal H}}\left(X\right) +\sum_{\ell =1}^{r}\sum_{i=2}^{r}\left\{ \left( B_{1i}^{{\cal H}^{r+\ell }}\right)^{2}-B_{11}^{{\cal H}^{r+\ell }}B_{ii}^{{\cal H}^{r+\ell
			}}\right\}. 
		\end{equation*}
		Since $\sum \limits_{\ell =1}^{r}\sum \limits_{i=2}^{r}\left( B_{1i}^{{\cal H}^{r+\ell }}\right)^{2}\geq \sum \limits_{i=2}^{r}\left( B_{1i}^{{\cal H}^{r+1}}\right)^{2}+\sum \limits_{\ell=2}^{r}\left( B_{1\ell }^{{\cal H}^{r+\ell }}\right)^{2}$, the aforementioned equation can be written as,
		\begin{equation}\label{eq-IC(4.7)}
			{\rm Ric}^{{\cal H}}\left( X\right) -{\rm Ric}^{{\cal R}}\left( F_{\ast}X\right) \leq \sum_{\ell =1}^{r}\sum_{i=2}^{r}B_{11}^{{\cal H}^{r+\ell}}B_{ii}^{{\cal H}^{r+\ell }}-\sum_{i=2}^{r}\left( B_{1i}^{{\cal H}^{r+1}}\right)^{2}-\sum_{\ell =2}^{r}\left( B_{1\ell }^{{\cal H}^{r+\ell}}\right)^{2}.
		\end{equation}
		Using (\ref{eq-IC(4.6)}) in (\ref{eq-IC(4.7)}), we derive
		\begin{equation*}
			{\rm Ric}^{{\cal H}}\left( X\right) -{\rm Ric}^{{\cal R}}\left( F_{\ast}X\right) \leq \sum_{\ell=1}^{r}\sum_{i=2}^{r}B_{11}^{{\cal H}^{r+\ell}}B_{ii}^{{\cal H}^{r+\ell }}-\sum_{i=2}^{r}\left( B_{1i}^{{\cal H}^{r+1}}\right)^{2}-\sum_{\ell =2}^{r}\left( B_{\ell \ell }^{{\cal H}^{r+1}}\right)^{2}.
		\end{equation*}
		Then defining $r$ real functions and performing computations similar to \cite[Theorem 4.9]{LLSV_2022}, we obtain the required inequality $(\ref{eq-IC(4.5)})$. We also note that the equality cases follow by the same computation as in \cite[Theorem 4.9]{LLSV_2022}.    
	\end{proof}
	
	\subsection{Example}
	Here, we give an example of a Riemannian map that satisfies the equality of Chen-Ricci and improved Chen-Ricci inequalities derived in Theorems \ref{main_thm_CRI} and \ref{main_thm_ICRI}, as follows.
	
	\begin{example}
		Let $$\left(M_1 = \left\{ (x_{1}, x_{2}, x_{3}, x_{4}) \in \mathbb{R}^{4} : x_{i} \neq 0~ \forall~ 1\leq i\leq 4 \right\}, g_{1} = e^{2x_1} \sum \limits_{i=1}^{2}dx_{i}^{2} + \sum \limits_{i=3}^{4} dx_{i}^{2}\right)$$ and $\left(M_2 = \left\{ (y_{1}, y_{2}, y_{3}, y_{4}) \in \mathbb{R}^{4}\right\}, g_{2} = e^{2y_1} \sum \limits_{j=1}^{2}dx_{j}^{2} + \sum \limits_{j=3}^{4} dx_{j}^{2}\right)$ be two Riemannian manifolds. Then by \cite{Polat_Meena}, $F: (M_1, g_1) \to (M_2, g_2)$ defined by
		\begin{equation*}
			F(x_1, x_2, x_3, x_4) = (x_1, x_2, x_3, 0)
		\end{equation*}
		is a Riemannian map with
		\begin{equation*}
			(\ker F_\ast)^\perp = \operatorname{span} \left\{h_i = \frac{\partial}{\partial x_i}\right\}_{i=1}^{3}, \quad (\ker F_\ast) = \operatorname{span} \left\{h_4 = \frac{\partial}{\partial x_4}\right\},
		\end{equation*}
		\begin{equation*}
			({\rm range}~ F_\ast) = \operatorname{span} \left\{ F_\ast h_j = \frac{\partial}{\partial y_j}\right\}_{j=1}^{3}, \quad ({\rm range}~ F_\ast)^\perp = \operatorname{span} \left\{\frac{\partial}{\partial y_4}\right\},
		\end{equation*}
		where $\left\{\frac{\partial}{\partial x_i} \right\}_{i=1}^{4}$ and $\left\{\frac{\partial}{\partial y_j} \right\}_{j=1}^{4}$ are bases of $T_p M_1$ and $T_{F(p)}M_2$, respectively. Further, we also have the only non-vanishing Christoffel symbols $\Gamma^{2}_{12}=\Gamma^2_{21}=1$, $\Gamma^{1}_{11} = 1$,
		$\Gamma^1_{22}= -1$ for $g_{2}$. Then by $(\ref{2ndff})$, we obtain
		\begin{equation*}
			(\nabla F_\ast) (h_i, h_j) = 0 ~ \text{for all $i,j \in \{1,2,3\}$}.
		\end{equation*}
		Hence $B^\mathcal{H}=0$. Subsequently, for $F$ the inequalities mentioned in $(\ref{eq-(3)})$ and $(\ref{eq-IC(4.5)})$ attain equality.
	\end{example}
	
	\subsection{Applications}\label{appl_part_CRI}
	In this subsection, we obtain the Chen-Ricci and improved Chen-Ricci inequalities for Riemannian maps to generalized complex and generalized Sasakian space forms. Toward this, we use the general forms established in Theorems \ref{main_thm_CRI} and \ref{main_thm_ICRI}. Also, for any ${\cal Z}\in \Gamma (TM_2)$, we decompose 
	\begin{equation}\label{decompose_GCSF_RM}
		J{\cal Z} ~\text{or}~ \phi {\cal Z} = P {\cal Z}+Q{\cal Z}, 
	\end{equation}
	where $P {\cal Z}\in \Gamma ({\rm range~}F_\ast)$, $Q {\cal Z}\in \Gamma (
	{\rm range~} F_\ast)^{\perp }$ with
    \begin{equation*}
        \|P F_\ast h_1\|^{2} = \sum \limits_{j=2}^{r}\left(g_{2}(F_{\ast }h_{1}, P F_{\ast }h_{j})\right)^{2} = \sum \limits_{j=1}^{r}\left(g_{2}(F_{\ast }h_{1}, P F_{\ast }h_{j})\right)^{2}.
    \end{equation*}
	
	First, we obtain the Chen-Ricci inequalities as follows.
	
	\begin{theorem}\label{appl_thm_CRI}
		Suppose $F$ satisfies the hypotheses of Theorem $\ref{main_thm_CRI}$. Then the following hold.
		\begin{enumerate}[$(i)$]
			\item If $(M_2(c_1, c_2), g_2, J)$ is a generalized complex space form, then
			\begin{equation*}
				{\rm Ric}^{{\cal H}}(X) \leq T_{GC} = \left( r-1\right) c_{1} + 3 c_{2} \|P F_\ast X\|^2 + \frac{1}{4}\left\Vert {\rm trace~} B^{{\cal H}}\right\Vert^{2}.
			\end{equation*}
			
			\item If $\left({M_{2}} \left( c_{1},c_{2},c_{3}\right), \phi, \xi, \eta, g_{2}\right)$ is a generalized Sasakian space form, then
			\begin{align*}
				&{\rm Ric}^{{\cal H}}(X) \leq \left\{\begin{array}{ll}
					& T_{GC}, ~{\rm if}~ \xi \in \Gamma({\rm range~}F_\ast)^\perp; \\& T_{GC}-\left( 1+T_\eta\right) c_3, ~{\rm if}~ \xi \in \Gamma({\rm range~}F_\ast),
				\end{array} \right.
			\end{align*}
			where $T_\eta = \left( r-2\right) \eta\left( F_{\ast }X\right)^{2}$.
		\end{enumerate} 
		The equality conditions are identical to those of Theorem $\ref{main_thm_CRI}$.
	\end{theorem}
	
	\begin{proof}
		Using (\ref{eq-(1)}), (\ref{Ricci_Curvatures}), (\ref{curvature_for_gcsf}), and (\ref{decompose_GCSF_RM}) we obtain,
		\begin{eqnarray*}
			&&{\rm Ric}^{{\cal R}}(F_{\ast}h_{1}) =\sum\limits_{j=2}^{r}c_{1} \{g_{1}(h_{j},h_{j})g_{1}(h_{1},h_{1})-g_{1}(h_{1},h_{j})g_{1}(h_{j},h_{1})\} \\&&+\sum\limits_{j=2}^{r}c_{2}\left\{ g_{2}\left( F_{\ast }h_{1},PF_{\ast }h_{j}\right) g_{2}\left( PF_{\ast }h_{j},F_{\ast}h_{1}\right) -g_{2}\left( F_{\ast }h_{j},PF_{\ast }h_{j}\right) g_{2}\left( PF_{\ast }h_{1},F_{\ast }h_{1}\right) \right. \\&&\left. +2g_{2}\left( F_{\ast }h_{1},PF_{\ast }h_{j}\right) g_{2}\left( PF_{\ast }h_{j},F_{\ast }h_{1}\right) \right\}.
		\end{eqnarray*}
		Equivalently,
		\begin{equation}\label{eq-GCSF-(2)}
			{\rm Ric}^{{\cal R}}(F_{\ast }h_{1})=\left( r-1\right)
			c_{1}+3c_{2}\sum \limits_{j=2}^{r}\left(g_{2}(F_{\ast }h_{1}, P F_{\ast }h_{j})\right)^{2}.
		\end{equation}
		Then $(i)$ follows by substituting (\ref{eq-GCSF-(2)}) in (\ref{eq-(3)}) as we can choose $h_{1}=X \in \Gamma \left( \ker F_\ast \right)^\perp$.\\\\
		In addition, using (\ref{eq-(1)}), (\ref{Ricci_Curvatures}), (\ref{curvature_for_gssf}), and (\ref{decompose_GCSF_RM}) we obtain,
		\begin{eqnarray*}
			&&{\rm Ric}^{{\cal R}}(F_{\ast}h_{1}) =\sum\limits_{j=2}^{r}c_{1}\left\{ g_{1}\left( h_{j},h_{j}\right) g_{1}\left( h_{1},h_{1}\right) -g_{1}\left( h_{1},h_{j}\right) g_{1}\left(h_{j},h_{1}\right) \right\} \\&&+\sum\limits_{j=2}^{r}c_{2}\left\{ g_{2}\left( F_{\ast }h_{1},PF_{\ast }h_{j}\right) g_{2}\left( PF_{\ast }h_{j},F_{\ast}h_{1}\right) -g_{2}\left( F_{\ast}h_{j},PF_{\ast }h_{j}\right) g_{2}\left( PF_{\ast }h_{1},F_{\ast }h_{1}\right) \right. \\&&+\left. 2g_{2}\left( F_{\ast }h_{1},PF_{\ast }h_{j}\right) g_{2}\left( PF_{\ast }h_{j},F_{\ast }h_{1}\right) \right\} \\&&+\sum\limits_{j=2}^{r}c_{3}\left\{ \eta \left( F_{\ast }h_{1}\right) \eta \left( F_{\ast }h_{j}\right) g_{2}\left( F_{\ast }h_{j},F_{\ast}h_{1}\right) -\eta \left( F_{\ast }h_{j}\right) \eta \left( F_{\ast}h_{j}\right) g_{2}\left( F_{\ast }h_{1},F_{\ast }h_{1}\right) \right. \\&&\left. +g_{2}\left( F_{\ast }h_{1},F_{\ast }h_{j}\right) \eta \left(F_{\ast }h_{j}\right) \eta \left( F_{\ast }h_{1}\right) -g_{2}\left( F_{\ast}h_{j},F_{\ast }h_{j}\right) \eta \left( F_{\ast }h_{1}\right) \eta \left(F_{\ast }h_{1}\right) \right\}.
		\end{eqnarray*}
		Equivalently,
		\begin{align}\label{eq-GSSF-(2)}
			&{\rm Ric}^{{\cal R}}(F_{\ast }h_{1})=\nonumber\\& \left\{\begin{array}{ll}
				\left(r-1\right) c_{1}+3c_{2}\sum\limits_{j=2}^{r}\left(g_{2}(F_{\ast }h_{1}, P F_{\ast}h_{j})\right)^{2}-\left( 1+\left(r-2\right) \eta \left( F_{\ast}h_{1}\right)^{2}\right) c_{3}, ~{\rm if}~ \xi \in \Gamma({\rm range~}F_\ast);\\
				\left(r-1\right) c_{1}+3c_{2}\sum\limits_{j=2}^{r}\left(g_{2}(F_{\ast }h_{1}, P F_{\ast}h_{j})\right)^{2}, ~{\rm if}~ \xi \in \Gamma({\rm range~}F_\ast)^\perp.
			\end{array} \right.
		\end{align}
		Then $(ii)$ follows by substituting (\ref{eq-GSSF-(2)}) in (\ref{eq-(3)}) as we can choose $h_{1}=X \in \Gamma \left( \ker F_\ast \right)^\perp$.
	\end{proof}
	
	Now, we obtain the improved Chen-Ricci inequalities as follows.
	
	\begin{theorem}\label{appl_thm_ICRI}
		Suppose $F$ satisfies the hypotheses of Theorem $\ref{main_thm_ICRI}$. Then the following hold.
		\begin{enumerate}[$(i)$]
			\item If $(M_2(c_1, c_2), g_2, J)$ is a generalized complex space form, then
			\begin{equation*}
				{\rm Ric}^{{\cal H}}\left( X\right) \leq IT_{GC} = \left( r-1\right) c_{1} + 3 c_{2} \|P F_\ast X\|^2+\frac{(r-1)}{4r} \left\Vert {\rm trace~}B^{{\cal H}}\right\Vert^{2}.
			\end{equation*}
			
			\item If $\left({M_{2}} \left( c_{1},c_{2},c_{3}\right), \phi, \xi, \eta, g_{2}\right)$ is a generalized Sasakian space form, then
			\begin{align*}
				&{\rm Ric}^{{\cal H}}(X) \leq \left\{\begin{array}{ll}
					& IT_{GC}, ~{\rm if}~ \xi \in \Gamma({\rm range~}F_\ast)^\perp; \\& IT_{GC}-\left( 1+T_\eta \right) c_3, ~{\rm if}~ \xi \in \Gamma({\rm range~}F_\ast),
				\end{array} \right.
			\end{align*}
			where $T_\eta = \left( r-2\right) \eta\left( F_{\ast }X\right)^{2}$.
		\end{enumerate} 
		The equality conditions are identical to those of Theorem $\ref{main_thm_ICRI}$.
	\end{theorem}
	
	\begin{proof}
		The proof follows by using (\ref{eq-GCSF-(2)}) and (\ref{eq-GSSF-(2)}) in (\ref{eq-IC(4.5)}).    
	\end{proof}
	
	\section{Conclusion with more applications and validations}\label{sec_conclusion}
	In Theorems $\ref{Theorem-GCRV}$,  $\ref{main_thm_CRI}$ and  $\ref{main_thm_ICRI}$, we have obtained general forms of the Chen-Ricci inequalities for Riemannian submersions and Riemannian maps. Further, by proving Theorems $\ref{appl_for_CRI_RS}$, $\ref{appl_thm_CRI}$ and $\ref{appl_thm_ICRI}$, we showed that these general forms yield new techniques that are easy, elegant, and fruitful for obtaining Chen-Ricci inequalities with other space forms. In view of $\cite{Olszak_1989}$, the item $(i)$ of Theorem $\ref{appl_for_CRI_RS}$ makes sense for non-constants $c_1$ and $c_2$ only when $m_1 = 4$. Similarly, the item $(i)$ of Theorems $\ref{appl_thm_CRI}$ and $\ref{appl_thm_ICRI}$ makes sense for non-constants $c_1$ and $c_2$ only when $\dim(M_2) = m_2 = 4$. However, for constants $c_1$ and $c_2$, the corresponding values of $m_1$ or $m_2$ could be greater than $4$. We note that one can directly obtain Chen-Ricci inequalities with several special space forms by substituting the values of $c_1$, $c_2$, and $c_3$ from Table \ref{table_3}. Hence, Theorems $\ref{appl_for_CRI_RS}$, $\ref{appl_thm_CRI}$ and $\ref{appl_thm_ICRI}$ are helpful to explore more particular cases, such as real, complex, real K\"ahler, Sasakian, Kenmotsu, cosymplectic, almost $C(\alpha)$ space forms. \\

	From \cite{Sahin_2013_I, Sahin_2010, Sahin_2013_SI, Sahin_2011, Park_Prasad, Tastan_Sahin_Yanan}, we have the notions of invariant, anti-invariant, Lagrangian, semi-invariant, slant, semi-slant, and hemi-slant Riemannian submersions from almost Hermitian manifolds onto Riemannian manifolds. In addition, from \cite{Lee_2013, Akyol_Sari_Aksoy, Erken_Murathan, Akyol_Sari, Sari_Akyol, Kumar_Prasad_Verma}, we have the notions of such structured Riemannian submersions from almost contact metric manifolds onto Riemannian manifolds. Similar to \cite[Chapter 6]{Sahin_book}, we can also extend the notions of such structured Riemannian maps to a complex or generalized Sasakian space form. These notions impose different-different conditions on the decomposition of the distributions $(\ker F_\ast = {\cal D}_1 \oplus {\cal D}_2)$ and $({\rm range~} F_\ast = {\cal D}_1 \oplus {\cal D}_2)$ with the complex structure $J$ or contact structure $\phi$. Using those conditions directly and putting the appropriate values of $c_1$, $c_2$ and $c_3$, one can also find Chen-Ricci inequalities for such structured Riemannian submersions and Riemannian maps. We observe that this approach successfully verifies the following existing inequalities.\\
	\begin{center}
		\begin{tabular}{cl}
			\hline
			\textbf{Reference(s)} & \textbf{Inequalities}\\
			\hline
			\cite{Akyol_Demir_Poyraz_Vilcu} & Equations (37), (38), and (62)\\
			
			\cite{Akyol_Poyraz} & Equations (44), (53), (66), (73), and (74)\\
			
			\cite{Aquib_2025} & Theorems 6, 7, 8, 9, and 10\\

			\cite{Aytimur_Ozgur, Aytimur_2023} & Theorems 4.1, 4.2, 4.3, 4.4, 4.5, and 4.6\\
			
			\cite{Gunduzalp_Polat_Filomat} & Theorem 4.1\\

			\cite{Gunduzalp_Polat_MMN} & Theorems 4 and 5\\
			
			\cite{LLSV_2022} & {\rm Theorems~ 3.1, 4.2, and 4.3}\\
			
			\cite{Polat_2024} & Theorems 4.1, 4.2, 4.3, 4.6, 4.7, and 4.8\\
			
			\cite{Poyraz_Akyol} & Theorems 3.4, 3.5, and 3.6\\
			\hline
		\end{tabular}
	\end{center}
	
	\vspace{0.2cm}
	
	Thus, the general forms obtained in this paper are powerful tools towards developing and verifying Chen-Ricci inequalities for Riemannian submersions and Riemannian maps. In the future, they may be explored for other remaining space forms.
					
					\section*{Acknowledgment}
					The second author is grateful to Professor Bayram \c{S}ahin for his comments that improved the initial version.
					
					\section*{Declarations}
					
					\begin{itemize}
						\item Funding: The author, R. Singh, was supported by the Human Resource Development Group (HRDG), Council of Scientific and Industrial Research (CSIR), New Delhi, India (File No: 09/0013(16054)/2022-EMR-I).
						
						\item Conflict of interest/Competing interests: Not applicable.
						
						\item Ethics approval: The submitted work is original and has not been submitted to more than one journal for simultaneous consideration.
						
						\item Consent for publication and participation: Not applicable.
						
						\item Availability of data, codes, and materials: Not available.
						
						\item Authors' contributions: All authors have equal contributions.
					\end{itemize}

	\noindent R. Singh\\
	Department of Mathematics, Banaras Hindu University, Varanasi, Uttar Pradesh-221005, India.\\
	E-mail: khandelrs@bhu.ac.in;	ORCID: 0009-0009-1270-3831\\
	
	\noindent K. Meena\\
	Department of Mathematics, Indian Institute of Technology Jodhpur, Rajasthan-342030, India.\\ E-mail: kirankapishmeena@gmail.com; ORCID: 0000-0002-6959-5853\\
	
	\noindent K. C. Meena\\
	Scientific Analysis Group, Defence Research and Development Organisation, Delhi-110054, India.\\
	E-mail: meenakapishchand@gmail.com; ORCID: 0000-0003-0182-8822\\

\addcontentsline{toc}{section}{References}
	\begin{thebibliography}{99}
	\bibitem{Akyol_Demir_Poyraz_Vilcu} Akyol, M.A., Demir, R., Poyraz, N.\"{O}., V\^ilcu, G.E.: Optimal inequalities for hemi-slant Riemannian submersions. Math., {\bf 10}(21) (2022), 1-18.
	
	\bibitem{Akyol_Poyraz} Akyol, M.A., Poyraz, N.\"{O}.: Sharp inequalities involving the Chen-Ricci inequality for slant Riemannian submersions. Bull. Korean Math. Soc., {\bf 60}(5) (2023), 1155-1179.
	
	
	\bibitem{Akyol_Sari_Aksoy} Akyol, M.A., Sar\i, R., Aksoy, E.: Semi-invariant $\xi^\perp$-Riemannian submersions from almost contact manifolds. Int. J. Geom. Methods Mod. Phys., \textbf{14}(5) (2017), 1-17.
	
	\bibitem{Akyol_Sari} Akyol, M.A., Sar\i, R.: On semi-slant $\xi^\perp$-Riemannian submersions. Mediterr. J. Math., \textbf{14}(6) (2017), 1-20.
	
	\bibitem{Alegre_Blair} Alegre, P., Blair, D.E., Carriazo, A.: Generalized Sasakian-space-forms. Israel J. Math., \textbf{141} (2004), 157-183.
	
	\bibitem{Aquib_Aldayel_Iqbal_Khan} Aquib, M., Al-Dayel, I., Iqbal, M., Khan, M.A.: Geometric inequalities and equality conditions for slant submersions in Kenmotsu space forms. AIMS Math., {\bf 10}(4) (2025), 8873-8890.
	
	\bibitem{Aquib_2025} Aquib, M.: Slant submersions in generalized Sasakian space forms and some optimal inequalities. Axioms, {\bf 14}(6) (2025), 1-15.
	
	
	\bibitem{Aytimur_Ozgur} Aytimur, H., \"{O}zg\"{u}r, C.: Sharp inequalities for anti-invariant Riemannian submersions from Sasakian space forms. {J. Geom. Phys.}, {\bf 166} (2021), 1-12.
	
	\bibitem{Aytimur_2023} Aytimur, H.: Curvature invariants for anti-invariant Riemannian submersions from cosymplectic space forms. Mediterr. J. Math., \textbf{20}(1) (2023), 1-17.
	
	\bibitem{Blair_2010} Blair, D.E.: {Riemannian Geometry of Contact and Symplectic Manifolds}. Birkh\"auser, Boston, (2010).
	
	
	
	\bibitem{Chen_2011} Chen, B.Y.: Pseudo-Riemannian Geometry, $\delta$-invariants and Applications. World Scientific, (2011).
	
	\bibitem{Chen_Blaga} Chen, B.Y., Blaga, A.M.: Recent developments on Chen-Ricci inequalities in differential geometry. Geometry of Submanifolds and Applications, Springer, (2024), 1-61.
	
	
	
	\bibitem{Deng_2009} Deng, S.: An improved Chen-Ricci inequality. Int. Electron. J. Geom., {\bf 2}(2) (2009), 39-45.
	
	\bibitem{Erken_Murathan} Erken, I.K., Murathan, C.: Slant Riemannian submersions from Sasakian manifolds. Arab J. Math. Sci., \textbf{22}(2) (2016), 250-264.
	
	\bibitem{Falcitelli_2004} Falcitelli, M., Ianus, S., Pastore, A.M.: {Riemannian Submersions and Related Topics}. World Scientific, River Edge, NJ, (2004).
	
	\bibitem{Fischer_1992} Fischer, A.E.: {Riemannian maps between Riemannian manifolds}. Contemp. Math., \textbf{132} (1992), 331-366.
	
	
	\bibitem{Gulbahar_Meric_Kilic} G\"{u}lbahar, M., Meri\c{c}, \c{S}.E., Kili\c{c}, E.: Sharp inequalities involving the Ricci curvature for Riemannian submersions. Kragujevac J. Math., {\bf 41}(2) (2017), 279-293.
	
	\bibitem{Gunduzalp_Polat_Filomat} G\"{u}nd\"{u}zalp, Y., Polat, M.: Chen-Ricci inequalities in slant submersions for complex space forms. Filomat, {\bf 36}(16) (2022), 5449-5462.
	
	\bibitem{Gunduzalp_Polat_MMN} G\"{u}nd\"{u}zalp, Y., Polat, M.: Some inequalities of anti-invariant Riemannian submersions in complex space forms. Miskolc Math. Notes, {\bf 23}(2) (2022), 703-714.
	
	
	
	
	
	\bibitem{Kumar_Prasad_Verma} Kumar, S., Prasad, R., Verma, S.K.: Hemi-slant Riemannian submersions from cosymplectic manifolds. Euro-Tbilisi Math. J., \textbf{15}(4) (2022), 11-27.
	
	\bibitem{Lee_2013} Lee, J.W.: Anti-invariant $\xi^\perp$-Riemannian submersions from almost contact manifolds. Hacett. J. Math. Stat., \textbf{42}(2) (2013), 231-241.
	
	\bibitem{LLSV_2022} Lee, J.W., Lee, C.W., \c{S}ahin, B., V\^ilcu, G.E.: Chen-Ricci inequalities for Riemannian maps and their applications. Contemp. Math., {\bf 777} (2022), 137-152.
	
	
	
	
	
	
	
	
	
	
	
	
	
	
	
	\bibitem{Neill_1996} O'Neill, B.: The fundamental equations of a submersion. Michigan Math. J., {\bf 13}(4) (1996), 459-469.
	
	\bibitem{Nore_1986} Nore, T.: {Second fundamental form of a map}. Ann. Mat. Pura Appl., \textbf{146} (1986), 281-310.
	
	\bibitem{Olszak_1989} Olszak, Z.: On the existence of generalized complex space forms. Israel J. Math., \textbf{65}(2) (1989), 214-218.
	
	\bibitem{Oprea_2010} Oprea, T.: Ricci curvature of Lagrangian submanifolds in complex space forms. Math. Inequal. Appl., {\bf 13}(4) (2010), 851-858.
	
	
	\bibitem{Park_Prasad} Park, K.-S., Prasad, R.: Semi-slant submersions. Bull. Korean Math. Soc., \textbf{50}(3) (2013), 951-962.
	
	\bibitem{Polat_2024} Polat, M.: B. Y. Chen-Ricci inequalities for anti-invariant Riemannian submersions in Kenmotsu space forms. Arabian J. Math., {\bf 13}(1) (2024), 181-196.
	
	\bibitem{Polat_Meena} Polat, M., Meena, K.: Clairaut semi-invariant Riemannian maps to K\"ahler manifolds. Mediterr. J. Math., \textbf{21}(3) (2024), 1-19.
	
	\bibitem{Poyraz_Akyol} Poyraz, N.\"{O}., Akyol, M.A.: Chen inequalities for slant Riemannian submersions from cosymplectic space forms. Filomat, {\bf 37}(11) (2023), 3615-3629.
	
	
	\bibitem{Sahin_2010} \c{S}ahin, B.: {Invariant and anti-invariant Riemannian maps to K\"ahler manifolds}. Int. J. Geom. Methods Mod. Phys., \textbf{7}(3) (2010), 337-355.
	
	\bibitem{Sahin_book} \c{S}ahin, B.: {Riemannian Submersions, Riemannian Maps in Hermitian Geometry, and Their Applications}. Elsevier, Academic Press, (2017).
	
	\bibitem{Sahin_2013_I} \c{S}ahin, B.: Riemannian submersions from almost Hermitian manifolds. Taiwanese J. Math., \textbf{17}(2) (2013), 629-659.
	
	\bibitem{Sahin_2010} \c{S}ahin, B.: Anti-invariant Riemannian submersions from almost Hermitian manifolds. Cent. Eur. J. Math., \textbf{8}(3) (2010), 437-447.
	
	\bibitem{Sahin_2013_SI} \c{S}ahin, B.: Semi-invariant Riemannian submersions from almost Hermitian manifolds. Canad. Math. Bull., \textbf{56}(1) 2013, 173-183.
	
	\bibitem{Sahin_2011} \c{S}ahin, B.: Slant submersions from almost Hermitian manifolds. Bull. Math. Soc. Sci. Math. Roumanie, \textbf{54}(1) (2011), 93-105.
	
	
	\bibitem{Sari_Akyol} Sari, R., Akyol, M.A.: Hemi-slant Riemannian submersions in contact geometry. Filomat, \textbf{34}(11) (2020), 3747-3758.
	
	
	
	\bibitem{Tastan_Sahin_Yanan} Ta\c{s}tan, H.M., \c{S}ahin, B., Yanan, \c{S}.: Hemi-slant submersions. Mediterr. J. Math., \textbf{13}(4) (2016), 2171-2184.
	
	\bibitem{Tricerri_1981} Tricerri, F., Vanhecke, L.: Curvature tensors on almost Hermitian manifolds. Trans. Amer. Math. Soc., \textbf{267}(2) (1981), 365-398.
	
	\bibitem{Tripathi_2011} Tripathi, M.M.: Improved Chen-Ricci inequality for curvature-like tensors and its applications. Differential Geom. Appl., {\bf 29}(5) (2011), 685-698.
	
	
	\bibitem{Vanhecke_1975} Vanhecke, L.: Almost Hermitian manifolds with $J$-invariant Riemannian curvature tensor. Rend. Semin. Mat. Univ. Politec. Torino, \textbf{34} (1975-76), 487-498.
	
	
	
	\bibitem{Yano_1984} Yano, K., Kon, M.: Structure on Manifolds. World Scientific, Singapore, (1984).
	
	
	
	\end{thebibliography}
\end{document}